\documentclass[11pt,notitlepage,twoside,a4paper]{amsart}
 \usepackage{amsfonts}
\usepackage[latin1]{inputenc}
\usepackage[cyr]{aeguill}
\usepackage{epigraph}

\usepackage{epsfig,fancyhdr,color}

\usepackage{xspace}
\usepackage[polutonikogreek,english]{babel}

\usepackage{amsmath,amssymb,enumerate}

\usepackage{epsfig,fancyhdr,color}

\usepackage{amssymb}
\usepackage{amsmath,amsthm}
\usepackage{latexsym}
\usepackage{amscd}
\usepackage{psfrag}
\usepackage{graphicx}
\usepackage[latin1]{inputenc}
\usepackage[all]{xy}
\usepackage[mathcal]{eucal}

\definecolor{NoteColor}{rgb}{1,0,0}

\renewcommand{\textsc}{\textcolor{red}}
\newcommand*{\tg}[1]{\textgreek{#1}}

\newtheorem*{theorem 1}{\rm\bf Proposition 1}
\newtheorem*{theorem 2}{\rm\bf Proposition 2}

\theoremstyle{definition}

\theoremstyle{remark}

\def\interieur#1{\mathord{\mathop{\kern 0pt #1}\limits^\circ}}

\title[René Thom: From mathematics to philosophy]{René Thom: From mathematics to philosophy}

\author{Athanase Papadopoulos}
\address{Athanase Papadopoulos,  Universit{é} de Strasbourg and CNRS,
7 rue René Descartes,
 67084 Strasbourg Cedex, France}
\email{athanase.papadopoulos@math.unistra.fr}
\makeindex
 
\date{\today}

\begin{document}

  \maketitle

  \epigraph{Toute ma métaphysique sous-jacente, c'est d'essayer de transformer le conceptuel en géométrique, le logique en dynamique.
  \\
(René Thom, from a letter to Claire Lejeune,  February 20 1980, quoted in \cite{Lejeune}).
}
  
  \begin{abstract}  In this chapter,
 I will discuss René Thom's approach to philosophy  based on his mathematical background. At the same time, I will highlight his connection with Aristotle, his criticism of the modern view of  science as a predictive process, his ideas on mathematical education, his position with respect to the French school of mathematics that was dominent in his time and his relationship with the philosophical community. I will also touch upon the connections between Thom's ideas and those of Leibniz, Riemann, Freud and others.
 
 The last version of this paper will appear as a chapter in the book Handbook of the History and Philosophy of Mathematical Practice (ed. Bharath Sriraman), Springer.

     \end{abstract}

    

  \section{Introduction}
  
   René Thom\index{Thom, René} belongs to a line of preeminent mathematicians who were thoroughly  engaged in philosophy.  After a brilliant but short career in  mathematics culminating in the Fields medal, he devoted himself to philosophy, in the context of a return to the ideas of Aristotle, the philosopher of nature, of whom he considered himself a heir.     
     Thom had a monumental program: to lay down the foundations of a new natural philosophy based on mathematics, and more particularly on geometry and topology. Natural philosophy, from his point of view, included physics, embryology, optics, hydrodynamics, biology, linguistics, ethology, semantics, psychoanalysis and other fields.

With his philosophical broad program, Thom was continuously an intriguing figure for linguists, psychoanalysts, biologists and  philosophers of sciences, let alone mathematicians. He opened up with these fields a dialogue which,  in the history of mathematics, is unprecedented by the width of its spectrum.  Philosophers of various sciences discussed and argued with him on technical grounds, defending their ways of thinking while trying to understand his original program, which, after all, seems quite natural for a mathematician: describing nature using the  fundamental ideas of topology, from the most basic notions of open and closed set, boundary, singularity of a differentiable function, to the more elaborate concepts of cobordism,  structural stability, analytic continuation, etc.

   Thom's philosophical program was initiated in his first published  book, \emph{Stabilité structurelle et morphogénèse} (Structural stability and morphogenesis) \cite{SS}, published in 1972, and it was continued in his \emph{Modèles mathématiques de la morphogenèse} (Mathematical models of morphogenesis) 
\cite{Morpho}, published in 1974 and his  \emph{Esquisse d'une sémiophysique: Physique aristotélicienne et théorie des catastrophes} (Smiophysics, a sketch; Aristotelian physics and catastrophe theory) \cite{Esquisse} published in 1988. He wrote other books with a strong philosophical flavor. He also 
 published a large number of papers on philosophy. I will quote several of them in the present chapter, and there are many more.

 Thom explained at several occasions the reasons for his passage to philosophy, and I shall dwell on them in what follows.\index{Thom, René}     He was also thoroughly involved in the reflection on mathematical education, and he participated to the major educational debates of his time, in particular, on the introduction of the so-called modern maths\index{modern maths} in schools, in France and elsewhere, and on the new university programs. His article  \emph{Les mathématiques modernes : une erreur pédagogique et
philosophique ?}  (Modern mathematics: an educational or philosophical error?) \cite{Erreur} (1970) was translated into several languages and published in various journals.
In both domains of pre-university and university programs, he was a defender of classical geometry, against the exaggerated weight given to algebra, set theory, information science and applied mathematics, and I will also comment on this in what follows.

  The plan of the rest of this chapter is the following. 
  In \S \ref{s:Vita}, I report on a few important moments from Thom's life. In particular, I review some events that were foreshadows of his future career as a mathematicians, and on certain factors that led him to abandon mathematical research and devote himself to philosophy. Section \ref{s:Bourbaki}, titled \emph{Thom, Bourbaki and mathematics education}, is intended to situate Thom with respect to the  Bourbaki group who was leading mathematical research and setting the rules for mathematical education in France, starting in the 1950s. At the same time, I indicate some of Thom's ideas on education and on mathematical writing. Section \ref{s:qualitative}, titled   \emph{Qualitative and quantitative}, is an overview of Thom's belief that, contrary to a common modern dogma, science is not quantitative but always qualitative, and of the debate that this belief generated among scientists.
   In \S \ref{s:rediscovers}, titled  \emph{Thom's rediscovery of Aristotle}, I comment on the relation of Thom's ideas with those of Aristotle, a filiation he has declared on many occasions and which he explains in several articles. In 
\S \ref{s:topology}, titled  \emph{Topology from an Aristotelian perspective}, I comment on the different assertions of Thom according to which he found in the philosophical ideas of Aristotle his own ideas on form, morphogenesis, and also the first notions of topology. 
In the next section, \S \ref{s:biology}, titled \emph{Aristotelian biology}, I explain how Thom\index{Thom, René} shared Aristotle's ideas on biology as a qualitative science, his ideas on hylemorphism,\index{hylemorphism} on homeomerous\index{homeomer} and anhomeomerous\index{anhomeomer} parts, and how he considered his program in this field to be also in the lineage of the one of the Stagirite.\index{Aristotle} Thom established relations between the mathematical notion of stratification\index{stratification} and other topological notions with ideas and concepts that Aristotle introduced in his work on the parts of the animals.
  Section
  \ref{s:singularities}, titled \emph{Philosophy of singularities}, is concerned with Thom's ideas on the role played by the mathematical theory of singularities in the general framework of prediction in science. The theory of singularities becomes, from this point of view, a philosophical topic. 
  Then follows a short section, 
\S \ref{s:ethology}, titled 
\emph{Topology and ethology}, which contains a few notes on Thom's ideas\index{Thom, René} on morphology\index{morphology} applied to this field.  The next two sections,
\ref{s:linguistics} and   \ref{s:psycho}   titled \emph{Linguistics and morphology} and  \emph{Psychoanalysis} respectively, are concerned with Thom's work on these two fields, which he viewed as domains of applications of his ideas on morphology.
Section
\ref{s:analogy},  \emph{A metaphysics of analogy},
 is a glimpse into Thom's ideas on analogy, again, in the tradition of Aristotle. In the last section  before the conclusion of this paper, titled \emph{A return to the classics},
I mention some ideas of  three mathematicians with ideas that were predecessors of Thom's reflections, namely, Leibniz, Riemann and Grassmann, pointing out relations with Thom's ideas.

   \section{\emph{Vita}: From mathematics to philosophy} \label{s:Vita}
  René\index{Thom, René!\emph{Vita}} Thom,\index{Thom, René} at several occasions, recounted some events  from his childhood that are related to the shift that occurred later in his intellectual activity, from mathematics to philosophy. I have already commented on this in the long biography \cite{AP-Biographie} that I devoted to him in the first chapter of the book  \cite{Thom-CNRS}. In this section, I will give a short \emph{Vita} of Thom, highlighting some episodes from his life that are connected with his journey from mathematics to philosophy.

              René Thom\index{Thom, René} was born in 1923 in Montbéliard, a small town in the East of France, now in the Burgundy-Franche-Comté region. The city was occupied by the Germans between 1940 and 1945,\footnote{The city of Montbéliard oscillated a few times between France and Germany. From 1042 to 1793, it was part of the Holy Roman Empire, after which it was attached to France. It  was  occupied by the Prussians during the Franco-Prussian War of 1870, and liberated one year later by the so-called Bourbaki army, named after Charles-Denis Bourbaki\index{Bourbaki, Charles-Denis} (1816-1897), a general of the army of Napoleon III and son of the Greek colonel Constantin-Denis Bourbaki.\index{Bourbaki, Constantin-Denis} The latter was educated in France and he served in the French army between 1787 and 1827.  The pseudonym Nicolas Bourbaki,\index{Bourbaki, Nicolas} which was chosen randomly, refers to the name of general Bourbaki.}  the period in which Thom\index{Thom, René} was finishing secondary school and entering higher education.  His parents owned a grocery store in Montbéliard. His father had a solid education, with a good knowledge of Latin, and he wrote poetry.   Thom used to say that he gained his open-mindedness from having met people with very different  backgrounds at the family grocery store.  His parents encouraged him and his siblings to get a higher education.

      Thom had the profile of a studious French schoolchild. From what he recounted, we know that he was not particularly brilliant at school, but some of his recollections from that period are harbingers of his life as a future mathematician, and I will mention  some of them.

      The first\index{Thom, René} reminiscence dates from age ten or eleven, when Thom had set himself the task of systematically exploring what every theorem he knew on the geometry of three-dimensional space becomes in dimension four. He writes about this in an interview published in 1989 \cite{Nimier}: 
   
      \begin{quote}\small
   [This] was, if I may say so, my first attempt to do something original; but it was a way for me to understand, say,  how a system of two planes in $\mathbb{R}^4$ was made,
etc., and I think I had arrived to a pretty good
intuition at that time. I was already seeing in four-dimensional space.\footnote{The translations from the French are mine. For some more subtle sentences, I have included the French original in a footnote.} 
\end{quote}

The second fact I mention took place when Thom was thirteen or  fourteen. We read in \cite{Porte} that he came across a book on differential calculus in the public library of his hometown and was seduced by it. But Thom's main interest was Euclidean geometry.\index{Euclidean geometry} He writes in \cite[p. 10]{Predire}:
``I had a certain taste for Euclidean geometry which immediately attracted me very much." On the other hand, he had no taste for algebra, which, he says, was akin to certain automatisms. He later maintained this opinion on this field, a field which he considered ``not very instructive", and he kept his ardor for Euclidean geometry,   ``our old Euclidean geometry", as he used to call it affectionately.  On the same subject, he writes in \cite{Thom-Logos}:

\begin{quote}\small
There was one encounter that was quite decisive, and this was with Euclidean geometry in the ninth and tenth grades,\footnote{I have transformed ``classe de troisième et de seconde" of the French educational system into their equivalent in the American one.} where I really developed a taste for geometric thinking and proof. I had made it a point of honor to find solutions to all conceivable problems in elementary geometry, in the geometry of the triangle, and this is why, since that time, I have kept a preference for Euclidean geometry that my colleagues and fellow mathematicians do not forgive me.
\end{quote}

 Besides geometry, Thom was intrigued by dynamics.\index{dynamics}\index{Thom, René} In the interview \cite{Nimier}, he recalls his early interest  and his motivation for this subject which, a few years later, became one of his favorite topics of reflection:

    \begin{quote}\small
    Around age seventeen, I started being interested in dynamics. I don't remember on what occasion it was, but I had handed in a paper to my mathematics teacher, where I talked about the eternal return seen from a dynamical point of view, the theories of the eternal return... It was the idea that one could have of a space-time, a universe in which there would be the eternal return, that is to say, where the dynamics would be periodic, but I believe that it was almost the first time that I really thought things in terms of dynamics.
\end{quote}

During his high-school years, Thom was equally fond of literature, poetry, history and science. He especially enjoyed reading Greek and Latin classical authors.  In 1939, the question arose as to the choice of the second part of the \emph{baccalauréat}.\footnote{In the old educational French system, the  ``baccalauréat", that is, the national high school diploma, had two parts. The ``baccalauréat 2ème partie" corresponded to the last year in high school.} Thom hesitated between the literature curriculum (which, in France, was called the \emph{philosophy} section) and the science one (called the \emph{elementary mathematics} section). He opted for the second. He explained his choice in his interviews with \'Emile Noël \cite [p. 9]{Predire}:  
\begin{quote}\small
We knew that this section offered more opportunities than the first one; it was perhaps an illusion, but we were convinced of that. Above all, we were in 1939, at the beginning of the war, and our parents, who had fought during the First World War, used to tell us:
``try be a artilleryman, you are less exposed than in the infantry! " To be an artilleryman, you had to have studied mathematics. This element probably weighed heavily in the birth of my mathematical vocation.
  \end{quote}
  
 Thom obtained his \emph{baccalauréat 2ème partie} in June 1940, in a \emph{lycée} in Besançon, the large city close to Montbéliard. He was sixteen.\footnote{Finishing school at age 16 is not unusual in France for a talented child.} The period was hard for the Thom family, as for many in Europe, because of the Second World War. Furthermore, in those years, studying at high school represented a major cost for families. Thom, as a gifted child from a modest background, had obtained a small government scholarship for his education. We read in an interview reported on by Michèle Porte \cite{Porte} that Thom's parents were really struggling to make their financial needs meet and  to put their children through education.\footnote{``Mon père a tiré le diable par la queue toute sa vie\ldots Mes parents se sont vraiment saignés pour me faire faire
des études."}

 After his \emph{baccalauréat 2ème partie}, Thom, together with his older brother, managed to cross the Swiss border, escaping from the Germans who had occupied Montbéliard and the region around it. Their stay in Switzerland was hard for they were penniless, and Thom recounts that they suffered there from starvation.  With the help of the French authorities, the two brothers returned to France, and they ended up in Lyon,  which was still in the free zone.\footnote{The city of Lyon was occupied later by the Germans, in November 1942.} There, with the help of a family friend, who at the same time was a remote cousin, René Thom enrolled again a high school, since he had nothing better to do, and he presented for the second time the \emph{baccalauréat 2ème partie}, this time in the philosophy section.\index{Thom, René}
 
    After Thom graduated anew from high school, the  preparation of the competition for  entering the \'Ecole Normale Supérieure was a natural choice for a talented young man like him, interested in science, literature and philosophy. In France, such a preparation is offered by some selective high schools,\footnote{This is again peculiar to the French system, where the preparation of the entrance exams to the so-called \emph{grandes \'Ecoles} is done in specialized classes provided by some high schools,\index{Thom, René} after the usual secondary school program.} and Thom was admitted at the \emph{Lycée Saint-Louis}, a prestigious secondary institution in Paris.\footnote{The \emph{Lycée Saint-Louis}, before Thom enrolled there, had among its pupils André Weil, Laurent Schwartz, Gustave Choquet and several other preeminent mathematicians.}   The preparation lasts in principle two years, in classes called \emph{classe de mathématiques supérieures} and \emph{classe de mathématiques spéciales}. Thom spent three years at the \emph{Lycée Saint-Louis}, as his first attempt to enter the \'Ecole Normale was unsuccessful and he had to repeat the year of \emph{mathématiques spéciales}.

 At several occasions, Thom recounted that right after he was admitted to the  \'Ecole Normale, he announced to the school director that rather than studying mathematics, his desire was to study philosophy of science, precisely, ``in the path of Cavaillès and Lautman". The acting director at that time was the physicist Georges Bruhat, who  strongly advised Thom not to follow this path, saying that it would be more useful for him to concentrate his efforts on the preparation of the \emph{agrégation}\footnote{This is a French diploma allowing, in principle, to become a teacher in a high school, but many students who plan to teach at university also pass this exam.}  in mathematics, cf. \cite{Nimier}.

   At the end of this paper, I have included an  appendix on the two philosophers of science  Jean Cavaillès\index{Thom, René}\index{Cavaillès, Jean} and Albert\index{Lautman, Albert} Lautman.\index{Cavaillès, Jean}\index{Lautman, Albert}

 Thom's teacher at the \'Ecole Normale was Henri Cartan,\index{Cartan, Henri} who took him under his tutelage and who played later an important role for him, especially at the beginning of his career as a mathematician.

Thom finished his studies at the\index{Thom, René}  \emph{\'Ecole}  in 1946, after having passed the competitive examination of the agrégation.  On the recommendation of Cartan, he was appointed to the University of Strasbourg.\footnote{In those times, one could enter academia in France without having a doctorate.}

       Thom was expected to submit a PhD thesis in the years that followed his appointment, and his first task was to find an appropriate research subject.  He chose the field of algebraic topology. He explained later that it was easier for a beginner like him to work in a new discipline such as this one, where the ground to be explored was virgin and where there were still easily accessible open problems \cite{INA}.

  Thom's work on cobordism,\index{cobordism} which was to ripen in his mind in the few years following his arrival in Strasbourg, inaugurated the implementation of the methods of differential topology in algebraic topology, so that the subject of topology of manifolds became a combination of the two, differential and algebraic topology. In fact, at the time Thom entered the field of topology, both points of view (the algebraic and the differential) were in a state of rapid expansion and were undergoing major transformations.\footnote{Thom entered topology just after a major revival of this field that started in the 1930s.  Topologists tried to classify manifolds and  were looking for algebraic invariants. The notions of fibre space\index{fiber space} and fibre bundle\index{fiber bundle}\index{Thom, René} became central in\index{algebraic topology} algebraic topology around the year 1950, and the topology of manifolds passed important milestones in the space of very few years: In 1956, Milnor\index{Milnor, John} proved the existence of exotic differentiable structures on 7-dimensional spheres.
In 1961, he disproved the so-called \emph{Hauptvermutung der kombinatorischen
Topologie} (``main conjecture of combinatorial
topology"), a conjecture formulated in 1908 by Steinitz et Tietze, asking whether  any two triangulations of homeomorphic spaces are  isomorphic after subdivision. The higher-dimensional analogue of the Jordan curve theorem (the so-called Jordan--Schoenflies theorem) was obtained in 1959-1960 (works of  Mazur,\index{Mazur, Barry} Morse\index{Morse, Marston} and Brown\index{Brown, Morton}). In 1961, Smale\index{Smale, Steven} proved the Poincaré conjecture for dimensions $\geq 7$. The same year, Zeeman obtained the corresponding result for dimension 5, and one year later, Stallings\index{Stallings, John} for dimension 6. Thom considers that he arrived in topology at the right time. He writes in  \cite{Thom-Logos}: ``I think that researchers who started in 
mathematics after 1960 have had infinitely more merit than
than people like me who came along at just the right time." }
    Algebraic topology and the study of differentiable functions occupied Thom\index{Thom, René}\index{Thom, René!doctoral dissertation} during his whole  mathematical career. He obtained his doctoral degree in 1951, with a dissertation\index{Steenrod squares} titled \emph{Espaces fibr\'es en sph\`eres et carr\'es de Steenrod} \cite{Thom-these}.
   Cartan,\index{Cartan, Henri} who was officially Thom's thesis advisor,\footnote{In those days, the task of the thesis advisor was generally that of a guarantor at the moment of the defense. In general, suggesting the subject of the thesis was not part of his role. At any rate, in Thom's case, it was not Cartan who proposed the topic of the work. I discussed this at length in the article \cite{AP-Biographie}.} became involved in this work after Thom\index{Thom, René} had written a first version, in order to correct the writing and make it intelligible. Several letters exchanged between Thom and Cartan, edited by André Haefliger\index{Haefliger, André} in the form of an article in volume I of Thom's Mathematical Collected Works \cite{Thom-oeuvres}, concern the finalization of the writing of this thesis. The article contains excerpts from these letters, which show that Cartan strongly committed himself during the last months before the defense, to make the thesis understandable. 
The general impression that emerges from this correspondence, as well as from some letters dating from the previous year concerning a paper that Thom had submitted to Cartan\index{Cartan, Henri} for publication in the \emph{Comptes Rendus de l'Académie des Sciences} (and which was eventually published in the proceedings of a conference), is that Cartan did not consider Thom capable of writing complete proofs. The correspondence between Thom and Cartan also shows that the latter consulted Samuel Eilenberg\index{Eilenberg, Samuel} about Thom's dissertation and that Eilenberg began by expressing serious reservations about the validity of the proofs. In particular, in a letter to Thom, Cartan writes that Eilenberg ``believes that it is absolutely impossible to know whether the theorems in question are correct, and consequently he considers that it would be very imprudent to publish them as they are in a thesis, without a true proof."
We note incidentally that among the changes that Cartan advised Thom to make to his thesis, he recommended the elimination of a complete chapter which seemed to him very unclear and almost impossible to make readable in order to defend the thesis in due time,\footnote{This was in the summer of 1951. Cartan had arranged an invitation for Thom\index{Thom, René}  to Princeton for the following fall, and the thesis had to be defended before.} hoping that Thom could at least work properly on the rest; in fact, this was a chapter on cobordism,\index{cobordism} a notion that Thom had introduced, and which became later the ground of the work for which he was awarded the Fields medal.

   From an exterior point of view, the years Thom\index{Thom, René}  spent in Strasbourg (1947-1963) were those of his purely mathematical activity. In his interviews \emph{Prédire n'est pas expliquer}, talking about this period, he summarizes it as, in the first place, that of an abandonment of philosophy in favor of mathematics, and also a period ending with a passage to mathematics as a tool for the description of the real world. This is how mathematics is used in \emph{catastrophe theory}, a theory which he started to develop under the name \emph{model theory},\footnote{The sub-title of Thom's book \emph{Stabilité structurelle et morphogénèse} is \emph{Essai d'une théorie générale des modèles} (An attempt for a general theory of models).} whose goal was to\index{morphogenesis}\index{structural stability} provide a mathematical framework to ideas from all domains of knowledge, including the physical, social and biological ones.
 He writes \cite[p. 24]{Predire}: 
   \begin{quote}\small
     From 1950 on, I abandoned these philosophical preoccupations to devote myself to mathematics. This lasted until 1956-1957. Then I went through a phase of depression: progress in mathematics was made by others. They led to theories so algebraically complicated that I could not follow them anymore. I gave up. \ldots
     But one must still do something! So I started to look for possible applications of the mathematical theories I knew. That is how I oriented myself towards catastrophe theory.\index{catastrophe theory} 
     
     \end{quote}

    By the end of his stay in Strasbourg,\index{Thom, René} Thom came across Christopher Zeeman's\index{Zeeman, Christopher} article \emph{Topology of the brain}, an article which opened for him the door to new possibilities of representing the phenomena of life by mathematical models. He  writes in his introductory talk of the 1982 symposium  \emph{Logos and Catastrophe Theory} \cite{Thom-Logos}:  
    \begin{quote}\small
    The idea of biological modeling, I owe it to a certain extent to Christopher Zeeman\index{Zeeman, Christopher} who, in 61, had published an article called Topology of the brain, where he had shown the great theoretical possibilities of mathematical modeling of the most complex biological activities. I was very impressed by this article at the time.
    \end{quote}

At the Edinburg ICM (14-21 August 1958), Thom\index{Thom, René}
was awarded the Fields medal for his work on cobordism.\index{cobordism} Later, he described this problem as ``to know when two manifolds can be deformed into each other without having singularities at any moment in this deformation". He considered that he had partially solved this problem. Concerning the prize itself, he declared in the interview \emph{Paraboles et catastrophes} \cite[p. 23]{Paraboles}: ``It was a bit of a good fortune. I seized the right moment by being able to exploit the techniques that Cartan, Ehresmann, Serre and many others had taught me." He wrote in his autobiographical article \cite{Atiyah} that ultimately, this prize left him with a bitter taste, because the work done was too little, compared to what others did after him. ``In this sense'', he says, ``this prize represented for me a certain fragility that the future made even more visible. The consequence, if there was one, was only the invitation to the new Institut des Hautes \'Etudes Scientifiques." Indeed,\index{Thom, René} Thom was offered a research position at this  newly founded institute.\footnote{The institute was founded by Léon Motchane (1900-1990), a Russian-born immigrant who started with small jobs in Switzerland and then Germany and who later became a wealthy businessman in banking and insurances. Motchane\index{Motchane, Léon} liked mathematics and he started studying them relatively late in his life, encouraged by Paul Montel. At  age 54, he obtained a doctorate, with Gustave Choquet as supervisor, and he decided to dedicate an important part of his time and fortune to the foundation of  a mathematical institute. He managed to obtain funding from several large private companies. The institute was founded in 1958, and Thom was among the first professors hired. The Institute was  situated in Paris, in a two-room apartment in the Thiers Foundation, located in a private mansion in the 16th arrondissement, before moving four years later to Bures-sur-Yvette, in the domain called Bois-Marie, which became its permanent location. Today, the French government is the main source of funding for the IH\'ES. Motchane\index{Motchane, Léon} represented Grothendieck for his reception  of the Fields medal at the Moscow ICM (1966).}      In  \cite[p. 27]{Predire}, he writes:
               \begin{quote}\small 
          I left Strasbourg and came, answering the call of the founder of the {\sc IH\'ES}.
  I had more leisure time, I was less preoccupied with teaching and administrative tasks. My purely mathematical productivity seemed to be on the wane and I began to focus more on the periphery, i.e., possible applications.   

  In addition to optics, I wondered if there were some possible applications to biology. I finally understood where these exceptional singularities came from: they are linked to the fact that the trajectories of light rays are special, because they satisfy a variational principle, Fermat's principle;\index{Fermat principle}\index{caustic} they have special properties that make caustics ``catch" singularities more easily than they should. It is the transfer of this idea that led me to the theory of catastrophes, to the mathematical part that leads to it.
            
      \end{quote}   
      Thom\index{Thom, René} declared at several occasions that his interest in forms arose from his observation of caustics.\index{caustic}   
 
Thom later had another ``revelation", during a visit to the Poppelsdorfer museum in Bonn, when he saw plaster cast models of the various phases of the development of a frog embryo\index{embryology} in the process of gastrulation. In particular, he noticed that some lines in these representations were, from the mathematical point of view, singularities that were familiar to him and he realized that the whole process may be interpreted using topology and geometry. He gave an interpretation of this beautiful geometry in terms of the unfolding of a wave front in an appropriate space that projects onto the ordinary space and later, he described all the process of development of the frog as a sequence of successive unfoldings of a singularity.\index{singularity!unfolding} This was the starting point for the application of his ideas on singularity theory to embryology \cite[p. 111]{Predire}.

     Thom's stay in Strasbourg lasted from 1947 to 1963, with a few breaks for visits he made to the United States, each  for several weeks, and an interruption during the academic year 1953-1954, which he spent in Grenoble.

   Thom introduced the word \emph{catastrophe}\index{catastrophe theory} in his book \emph{Stabilit\'e structurelle et morphogén\`ese} (1972), to denote an abrupt phenomenon that appears and disappears continuously. The work quickly became the reference for the theory that became known under the name \emph{Catastrophe theory}.
  Speaking of this theory, Thom recalls the role of the optical experiments he conducted in 1959-1960 in Strasbourg as an impetus. He declares in  his interview \cite[p. 26]{Predire}:   
  \begin{quote}\small
  I came to [catastrophe theory] in a rather natural way: it was an evolution that led me, starting from a purely mathematical problem, the one called the generic singularities of a map, to see whether this theorem had physical applications. 
  I was then at the University of Strasbourg and I had the opportunity to do some optical experiments. A colleague, who was a physicist,  lent me some instruments, a spherical mirror, a prism and a dioptre. I was able to produce some caustics, and I made them vary by changing a little bit the position of the parameters; I observed how they were deformed. That is precisely the thing that interested me.  [\ldots] From there I started the theory of catastrophes: I realized that there were types, variations of caustics, that I did not foresee, and that I had to explain to myself the appearance of these singularities. It took me two or three years to understand where this came from. This is a banal phenomenon, but that is where I started from.
  \end{quote}

 Optical caustics are curves that appear as envelopes of a beam of light rays reflected in a curved mirror. Later, in  1991, in his  lecture titled \emph{Leaving mathematics for philosophy} \cite{Leaving}, Thom recalled that his first surprise during his experiments on optics was to see that the caustics obtained using a spherical mirror or a linear diopter give rise to singularities that he had not expected to see.

 In 1961, during one of his trips to the United States, Thom met Salomon Lefschetz,\index{Lefschetz, Salomon} who was working on structural stability,\index{structural stability} a notion that was introduced in the 1930s by  Aleksandr Andronov and Lev Pontryagin. 
    Several years later,  in his introductory talk of the 1982 conference \emph{Logos et théorie des catastrophes} (Logos and catastrophe theory) \cite{Thom-Logos},  Thom\index{Thom, René}  wrote, regarding Lefschetz : ``He understood very well the crucial importance of this notion and he did a lot to develop this theory in the Western circles." Later on, Thom used this concept 
 in a decisive way in his mathematical and philosophical work, and it became a central element of his theory of morphogenesis.  In fact, the title of his first book, published in 1972, whose aim is to provide a mathematical formalism  to tackle general morphogenesis problems,
 is precisely \emph{Stabilité structurelle et morphogénèse} \cite{SS}. He writes in the introduction: ``This book is the first systematic attempt to think in geometrical and topological terms the problems of biological regulation, as well as those posed by the structural stability of any form."\index{Thom, René}
 Thom, whose main objective became  to understand form in nature, considered that only structurally stable morphologies and phenomena can be understood. He even tried to formalize thermodynamics using the notion of structural stability, although he realized later on that he failed, and that some weaker notion of structural stability was needed.\footnote{Thom wrote:  ``initially, the theory of structural stability\index{structural stability} seemed to me to be of such scope and generality 
that with it,   I could hope to replace thermodynamics by geometry, 
to geometrize in a certain sense thermodynamics, eliminating
from thermodynamical considerations all the stochastic and measurable aspects, to retain only the corresponding geometric  
characterization of the attractors. It is certain that
the instability phenomena of the attractors that have been discovered since then
show that such a hope is false or, at least, that it would be necessary
 to weaken considerably the notion of structural stability. " \cite[p. 587]{Thom-Logos} (1989).} 
 
  Let me make a parenthetical remark, concerning physics regarded by philosophers from Greek Antiquity, since, as we shall recall below, Thom often referred to them. Plato considered that there cannot be any knowledge of nature, that is, that physics does not exist, since the physical world is permanently subject to change, and it is impossible to have a real knowledge of something which is perpetually changing. Aristotle,\index{Aristotle} on the contrary, in his \emph{Physics}, was puzzled by evolution, and by things that are changing. His \emph{Physics} is precisely a physics of motion and of change. Structurally stable phenomena makes a compromise between changing and  unchanging objects.  A dynamical system is structurally stable if, at each moment of its evolution, there is an interval of time around this moment such that  the system obtained at any given time  in this interval 
 is identical to the original system.\footnote{The mathematical formulation of this notion asserts the existence of a conjugacy between a system and close-enough systems so that in some sense the two systems are similar.} 
 
 Thom's\index{Thom, René} main philosophical program, that of describing nature and language in geometric terms, is already largely initiated in his book \emph{Stabilité structurelle et morphogénèse} (1972). He writes in the Preamble: ``This is a work that wishes to place itself in the line of a deceased discipline, namely \emph{Natural philosophy},"  and the reference to Aristotle is clear.
 A good half of this book is devoted to questions of embryology, where Thom proposes a formulation of this field using a geometric language.

In his second book, titled \emph{Modèles mathématiques de la morphogenèse} \cite{Morpho}, published in 1974, Thom defines morphogenesis as ``the creative or destructive process of forms"  \cite[p. 10]{Morpho}.  Thom declares that his research on morphogenesis is based on two sources: on the one hand, his research in differential topology on structural stability, and on the other hand, his reading of embryology treatises. He  considered that the applications of the theory he proposes go far beyond embryology and even biology, to include subjects such as physical chemistry, optics, hydrodynamics and others.

   At IH\'ES, Thom\index{Thom, René} continued publishing on mathematics, but not in a classical style, that is, not on new theorems and proofs. Most of his mathematical works were addressed to biologists, philosophers, linguists\index{linguistics} and other non-mathematicians, and some of these articles were published in philosophical journals. Thom declared later that his drift towards biology and the other sciences was due in part to the presence of  Alexander Grothendieck\index{Grothendieck, Alexander} at IH\'ES. He writes on this subject, in an autobiographical text \cite[p. 82]{Atiyah}: 
    \begin{quote}\small  
 [\ldots]   Relations with my colleague Grothendieck were less agreeable for me. His technical superiority was crushing. His seminar attracted the whole of Parisian mathematics, whereas I had nothing to offer. That made me leave the strictly mathematical world and tackle more general notions, like the theory of morphogenesis, a subject which interested me more and more and led me towards a very general form of ``philosophical" biology. The director, L. Motchane, had no objection (or if he did he kept it for himself). And so I ran what was known as a ``crazy" seminar, which lasted the best part of my first year at the institute. Three or four years passed, during which I returned to Aristotle and classical Greek.
        \end{quote}
        
%
%
%

          \section{Thom,  Bourbaki and mathematical education}\label{s:Bourbaki}
      
At the \'Ecole Normale Supérieure, Thom received from  Cartan an education  in the purest Bourbaki tradition. Soon after, in Strasbourg, he was seduced by ideas of Ehresmann, another eminent member of the Bourbaki group. But later, he became no less than an acerbic critic of the spirit that reigned in that group to which he owed his mathematical training. In his interviews \emph{Paraboles et catastrophes} \cite[p. 23]{Paraboles}, reporting on a  Bourbaki\index{Bourbaki, Nicolas} session\index{Thom, René} to which, at the request of Cartan,\index{Cartan, Henri} he participated while he was still  a student at the \'Ecole Normale, he writes: 
 \begin{quote}\small
 According to the tradition, I had to attend a presentation of the texts and I was asked for my opinion in case of disagreement.  \ldots  I was tired of listening (it was a bit boring!) and sometimes I even fell asleep during the sessions. In a sense, this was not a bad thing after all: later, in fact, I became an opponent of the Bourbaki view of mathematics. \ldots  They were really ultra-formalists! They presented mathematics in a rigorous,\index{rigor} ascetic manner; but in doing so, they ended up dealing only with those portions of mathematics that were appropriate to this rigorous presentation. For this reason, much of the mathematics that was practiced at the time could only escape their attention. \ldots 
Even today there are large sections of\index{Thom, René} mathematics---the calculus of variations, partial differential equations, qualitative dynamics, etc.---that are considered as not ``clean"  enough (or ``dead"!) to be included in the Bourbaki\index{Bourbaki, Nicolas} anthology. Thus, there is a genuine incompatibility between the Bourbaki repertoire and living mathematics. Bourbaki has embalmed the mathematics that he has reduced, so to speak, to a mummy!
\end{quote} 
 
 Rigor,\index{rigor} systematization, and axiomatization were among the notions at the foundation of the Bourbaki style.\index{Bourbaki, Nicolas}  In the article \emph{Modern mathematics: a pedagogical and philosophical error?}  \cite{Erreur}, Thom expresses a certain retreat with respect to formalization, and he speaks of the limits of axiomatization.  He protests against the disappearance of traditional Euclidean geometry (the geometry of the triangle)\index{Euclidean geometry} from French curricula, in favor of set-theoretic and algebraic theories, where students no longer take  initiatives to solve construction problems, but learn abstract notions and apply blind algebraic recipes from a pre-established scheme. He then comes to the notion of ``meaning", and he raises the question of why we do mathematics. Is it a gratuitous game, a random product of our brain activities, or is it a way to comprehend the universe? At the same time, he\index{Thom, René} wonders about the meaning of the word ``rigor" in mathematics.\index{rigor} His answer is: ``A proof is rigorous, if it arouses, in any sufficiently educated and prepared reader,  a state of evidence that leads to adhesion" (p. 230). This is not the orthodox view on rigor, where a proof is rigorous\index{rigor} if it consists of a logical sequence of implications, starting from well-known statements.
He recalls that the etymology of the word ``theorem " is ``vision", declaring that ``there is no need for great axiomatic constructions or refined conceptual machineries to judge the validity of a reasoning". 
He writes (p. 232): 
\begin{quote}\small
There has been a lot of talking in recent years about the importance of axiomatization as an instrument of systematization and discovery. Yes, an 
instrument of systematization, certainly; but for what concerns discovery,  the matter is more than
doubtful. It is characteristic that, from the immense effort of systematization of Nicolas
Bourbaki\index{Bourbaki, Nicolas} (which, by the way, is not a formalization, because Bourbaki uses a non-formalized metalanguage), no new theorem of any importance has emerged. And if researchers in 
mathematics refer to Bourbaki,\index{Bourbaki, Nicolas} they find their food much more often in the exercises---where the author has pushed the concrete material---than in the deductive body of the text. This must be said quite clearly: axiomatization is a research intended for specialists, and has no place either in secondary teaching nor at 
university (except, of course, for professionals who wish to specialize in the study of the foundations).  This is why accusations of inconsistency addressed to Euclidean geometry are in fact irrelevant at the level 
of the local intuitive validity of the reasoning, which is the only level that matters.
\end{quote} On p. 236 of the same article, he writes: 
\begin{quote}\small
It is not certain that even in pure mathematics, every deduction can have a set-theoretic model. The ill-tamed paradoxes that undermine formal set theory are there to remind the mathematician of the dangers of the exaggerated use of these seemingly innocent symbols. Perhaps, even in mathematics, quality remains, and resists any set-theoretic reduction. The old\index{Thom, René} Bourbakist\index{Bourbaki, Nicolas} hope, to see the mathematical structures  emerge naturally from the hierarchy of sets,  their subsets and their combinatorics, is probably only a fantasy. No one can reasonably escape the impression that the important mathematical structures (algebraic structures,
topological structures) appear as data fundamentally imposed by the external world, and that their irrational diversity finds its   justification only in reality.
\end{quote}

Thom often referred to Gödel\index{Gödel, Kurt} to remind us that absolute\index{rigor} rigor in mathematics does not exist, the latter having shown that a completely formalized system of arithmetic (like a machine) is either inconsistent (leading to a contradiction), or incomplete.\footnote{I would like to point out the interesting paper \cite{FN} by Farmaki and Negrepontis titled \emph{The paradoxical nature of mathematics}, in which the authors argue that the deductive strength in Mathematics is strongly related to its paradoxical nature, and in fact, that it comes, maybe exclusively, from its proximity to the\index{contradictory nature of mathematics} contradictory.}

 In a dialogue bearing the title \emph{Métaphysique extrême} (Extreme metaphysics), Thom\index{Thom, René} declares \cite[p. 86]{Onicar}:
   ``Formalists are people who always tell you: `Oh, natural language is horrible, it tolerates every kind of ambiguity, it is not possible to do mathematics with it.' I have never done anything but natural language in mathematics, plus a few symbols from time to time."

 In a series of lectures he gave in 1981, titled \emph{Mathématiques essentielles}\index{essential mathematics} \cite{Essentielles}, given at the Solignac Abbey, a medieval monastery in the Limousin (South of France), the main idea that Thom tried to communicate is that ``what justifies the `essential' character of a mathematical theory is its ability to provide us with a representation of reality."  He declares \cite[p. 2-3]{Essentielles} that the various branches that constitute mathematics are not all essential and that the possibility of abstracting mathematical entities from concrete situations comes from the fact that mathematics provides us with a representation of the real. He continues:
\begin{quote}\small
 By ``real" I mean both the reality of the external world---whether it is given to us by the immediate perception of the world around us, or by a mediate construction as is scientific vision; that is to say, that the ``real", in its most immediate form, is also apprehended by introspection---as an ``immediate data of consciousness".  
 \end{quote}

Thom reiterated this thesis in a talk he gave in 1982 at the \'Ecole Normale, titled  \emph{Les réels et le calcul différentiel, ou la mathématique essentielle} (The reals and differential calculus, or essential mathematics) \cite{essentiel} and in several other writings. He was often accused of giving non-rigorous demonstrations,\footnote{See e.g. \cite{JQ} and Thom's response in \cite{Collectif}.} but these non-rigorous demonstrations were part of his style, which was that of a continuous search for meaning instead of rigor, and a refusal of any conventional way of writing mathematics. 
 In his interview with Jacques Nimier \cite{Nimier}, he expresses his doubts and hesitations,  which were part of his philosophy:
  \begin{quote}\small
 There are clean theories and dirty theories, and personally I always have more sympathy for a dirty theory. The clean theories are the theories where things are well presented, where the concepts are clearly defined, the problems more or less well defined as well.\index{rigor} Whereas the dirty theories are the theories where one does not know very well where to go, one does not know how to organize things and where the main directions are, etc. From this point of view, in fact, I have never been a Bourbakist,\index{Bourbaki, Nicolas} because Bourbaki likes clean things; whereas I think that it is necessary to get one's hands dirty, and sometimes even dirtier in mathematics.
 \end{quote}

        As a passionate of topology in the sense of the theory of forms and their deformations, Thom was not attracted by algebra, nor by number theory, nor by set theory whose formalism had invaded  mathematical writing. He adopted a way of thinking that included mathematics in a very broad perspective, that of motion, form and change of form.\index{form} Gradually, he integrated the mathematical objects he had been working on into his new framework of thought, that of natural philosophy, where mathematics becomes the model of the physical and biological world.  In his article \emph{Les racines biologiques du symbolisme}  (The biological roots of symbolism) \cite{racines} (1978), he\index{Thom, René}  returns to the opposition between algebra and geometry, this time associating topology with geometry: 
 \begin{quote}\small
 There is a certain opposition between geometry and algebra. The fundamental material of geometry, of topology, is the geometrical continuum; pure and unstructured extent is \emph{par excellence} a ``mystical" notion.   Algebra, on the contrary, shows a fundamentally ``diaeretic" operative attitude.   Topologists are the children of the night; algebraists, on the other hand, wield the knife of rigor with perfect clarity.
\end{quote}

      Thom's judgment on algebra\index{rigor}\index{algebra} was shared by some other preeminent mathematicians. Michael Atiyah, referring to Goethe, considers that when a mathematician solves a question using algebra, he sells his soul. He writes in his article \emph{Mathematics in the 20th Century} \cite{Atiyah-math}:  
      
 \begin{quote}\small Algebra\index{Thom, René} is the offer made by the devil to the mathematician. The devil says: ``I will give you this powerful machine, it will answer any question you like. All you need to do is give me your soul: give up geometry and you will have this marvelous machine". [\ldots] When you pass over into algebraic calculation, you essentially stop thinking; you stop thinking geometrically, you stop thinking about the meaning.
\end{quote}

Mostly in reaction to the  proponents of ``modern mathematics"\index{modern maths} that imposed new programs in school and university curricula, Thom published, in 1962, an article on envelopes \cite{T1962}, a classical subject which originates in the work of the 17$^{\mathrm{th}}$ century mathematician Christiaan Huygens.\index{Huygens, Christiaan} The subject had been in the high school curriculum before the modern math reform, but had disappeared, certainly because of a difficulty it presented in comparison with subjects like linear algebra or analysis: in the study of envelopes, there is no general method, and one had to reason case by case. This was the kind of mathematics that Huygens preferred, and Thom followed him.
 Furthermore, this subject accounts for phenomena of life: Thom,\index{Thom, René} again like Huygens, had observed the envelopes\index{envelopes} of light rays undergoing reflection or refraction on a surface or curve and the structure of the singularities\index{singularity} they form. He was interested in the distinction between those that are stable or generic and those that are not, as well as in problems of classical plane geometry that the theory of envelopes leads to.\index{envelopes} Furthermore, Thom included the study of envelopes in that of generic singularities of differentiable functions, one of his favorite subjects. In fact, in his auto-biographical article \cite{Problemes}, Thom writes that he started the study of differentiable functions in order to understand envelopes.\index{envelopes}  
 
\section{Qualitative and quantitative}\index{quantitative!science}\index{qualitative!science}\label{s:qualitative}

   From\index{science!quantitative}\index{science!qualitative} the beginning of the 1970s, Thom\index{Thom, René}\index{analogy!René Thom} became involved in epistemology, in particular in the question of what is expected from science. In his books and articles, he expressed his ideas on the predictive attribute of theory, on the role of experimentation, and on analogy as a fundamental tool for reflection. From epistemology, he was 
gradually led to the great problems of philosophy. He turned more and more to Aristotle,\index{Aristotle} with his categories\index{Aristotle!categories} and genera,\index{Aristotle!genera} expressing reservations about mathematical logic as the basis of thought, stressing the fact that it may lead to contradictions and making clear his preference to Aristotelian logic, a logic used in a natural way, and based on reality. 

Thom liked Aristotle's\index{Aristotle} informal and lively style. He writes in his \emph {Esquisse d'une sémiophysique}: ``The very style of Aristotle,\index{Aristotle} far from the axiomatic\index{Thom, René} precision expected of the logician, is that of a thinking that seeks itself, made of returns on oneself, always in struggle with its object; this testimony of a constant effort of thinking filled me with a kind of immense sympathy" \cite[p. 13]{Esquisse}.
He emphasized at several occasions that the fact that the theory he was proposing had no quantitative character was an echo of Aristotelian biology. This is the subject of his article \emph{Qualitative and quantitative in evolutionary theory with some thoughts on Aristotelian Biology} \cite{Qualitative}. He used to describe Aristotle as ``the apostle of the qualitative against the quantitative" \cite[p. 172]{Esquisse}.
In his 1988 article \emph{Un rapport problématique: la rencontre théorie-expérience} (A problematic relationship: the theory-experiment encounter) \cite{Exp}, he contests the frequently repeated postulate that modern science was born at the beginning of the 17$^{\mathrm{th}}$ century as a result of the development of the experimental method. Science, he insists, is the consequence of theoretical reflection, and no experimentation is useful without a causal analysis that sets its limits and allows for an adequate interpretation of the results.

 Thom\index{Thom, René} declared in 1984, in a communication to the Académie des Sciences titled \emph{La méthode expérimentale: un mythe des épistémologues (et des savants) ?} (The experimental method: A myth of epistemologists (and Scientists)?) \cite{methode} that only two types of causal analyses are available for science:  
  the first one, of Aristotelian origin, is that of efficient causes,  based on natural language, explaining the phenomena by \emph{ad hoc}  entities, and the second one, linked to mathematics, essentially reducing the question of causality to the solution of a differential system whose solutions are completely determined by the initial solutions. In biology, he says, the first scheme is massively used to deduce precepts that have no demonstrative value. In the same article, Thom writes that  the expression ``experimental method" is not only a myth, but also an oxymoron: it contains a contradiction.  Indeed, he recalls that, while the word ``method" refers to a direction of thinking,\footnote{Etymologically, the word ``method" consists of two words:  the word is a combination of the words \tg{met'a}, which means ``after", or, ``to follow" and    \tg{<od'os}, which means a ``way", or a ``path".} the word experiment originates in the latin word \emph{experientia}, which means a trial.
 Thus, rather than a ``method", experimental\index{science!experimental} science is a \emph{pratice}, a practice, which according to Thom,\index{Thom, René} comes under the heading of exploration carried out at random, if it is performed without a theory supporting it. The aim of experimentation would be, according to him, to find ``that which one has not looked for", and, as such, it is illusory in many respects. He considers that approximation should  not be part of science as an enterprise aiming at constituting a knowledge with universal validity. Specifically, Thom criticizes what he calls ``contemporary biology," which uses the experimental method with no scruple.  
 
 A fierce debate followed Thom's lecture at the Academy, notably with the physicist Anatole Abragam, a fervent supporter of the experimental method, and who was also an academician. The text of Thom's lecture as well as Abragam's answer are published in the volume \emph{La philosophie des sciences aujourd'hui} (Philosophy of sciences today), edited by J. Hamburger \cite{Hamburger}.

Thom returned to this question several times in his writings. In his interviews \cite[p. 130]{Predire}, he writes:  ``If one reduces science to a set of recipes that work, then, intellectually speaking,  we are not  better than a rat who knows that when a lever is pressed, food will fall into his bowl. The pragmatist theory of science brings us back to the situation of the rat in its cage." This is illustrated on the book cover \cite{Predire}

%
%
 
In the same book \cite[p. 79]{Predire}, Thom\index{Thom, René} returns to a quote by Rutherford,  ``Qualitative is nothing but poor qualitative", which he strongly contests. He is convinced that the contrary is true.  He writes:  ``All of topology is to be added to the chapter of examples of this conviction. How is a sphere different from a ball? This is really not quantitative. How is a circle different from a disk? It is not a question of quantity, it is a problem of quality."
 
One characteristic of Thom's theory of morphogenesis\index{morphogenesis} which will be criticized by the proponents of a theory of pure experimentation, is that it is a completely geometrical theory, ``independent of the substrate of forms and of the nature of the forces that create them" \cite[p. 8]{SS}. According to him, the fact that a theory is independent of the substrate makes it applicable  at the same time to domains which are  a priori very different from each other, like biology and linguistics. The ground space is simply an open set in a topological space, i.e., an abstract substrate, devoid of any matter. Thom writes: ``This may seem difficult to admit, especially for experimentalists who are accustomed to cutting to the chase and continually struggling with a reality that resists them."

Thom said later that he ended up being a philosopher of science because of the need to respond to the criticisms directed towards catastrophe theory. He declared in an interview \cite[p. 90]{Predire}:
\begin{quote}\small
I was indeed subjected to criticisms that were on the epistemological  ground. I devoted myself for a while to the philosophy of science before engaging in more general philosophy [\ldots]. I started somehow from the validity of the theory of catastrophes, and I came to be interested in the position of science in general, in what can be expected from it from the point of view of knowledge. It was there that I began to develop critical positions towards the so-called experimental method and the general belief  in the virtues of experimentation as leading us to progress. I believe that experimentation by itself can hardly lead to progress. [\ldots] But I found myself hitting  a kind of wall: in a science like biology, for example, people refuse the theoretical necessity of imagination.
\end{quote}
 
     \section{Thom's rediscovery of Aristotle}  \label{s:rediscovers}

Thom\index{Thom, René} has always had a deep reverence for the philosophers of classical Greece.\index{Aristotle} From the end of the 1960s on, he constantly referred to them. In his book \emph{Structural stability and morphogenesis}  \cite{SS} which we already mentioned several times, he writes: ``In fact, all the fundamental intuitions in morphogenesis are already to be found in Heraclitus; my unique contribution will be to have placed them in a geometrical and dynamical framework which one day will make them accessible to quantitative analysis."  

Among the great figures of antiquity, Aristotle\index{Aristotle} occupies a central place in Thom's\index{Thom, René} work.  On several occasions, he expressed his surprise for having found in Aristotle's writings ideas that he had previously held himself. He also shared with the Philosopher a boundless curiosity and a desire to understand all the phenomena of nature. In his book, \emph{Semiophysics, a sketch},  which bears the subtitle \emph{Aristotelian physics and catastrophe theory}, he writes: ``It is only rather recently, almost incidentally, that I discovered Aristotle's work. I was almost immediately fascinated by this reading" \cite[p. 12]{Esquisse}.    In his \emph{Esquisse}, published in 1988, he writes \cite [p. 12]{Esquisse}:
\begin{quote}\small 
 [\ldots] If I add that I found in Aristotle\index{Aristotle}\index{Thom, René} the concept of genericity\index{genericity} (\tg{<ws >ep`i tò pol'u}), the idea of ``stratification'' that we can glimpse in Aristotle the biologist in the decomposition of the organism into homeomers\index{homeomer} and anhomeomers,\index{anhomeomer} the idea of the decomposition of the genus into species as an image of bifurcation, we will agree that there was something to be astonished about.
 \end{quote} 
 
 Thom also discovered in Aristotle the idea of significant fact and significant form, which played an important role in his philosophy.   He even found the essence of the mathematical notions of cobordism,\index{cobordism}  genericity\index{genericity} and stratification\index{stratification}  in the writings of the Stagirite,\index{Aristotle!biology} in particular in his biology.\index{biology!Aristotle}  I shall talk about this in more detail in the next sections. 
   The titles of some of Thom's articles are significant\index{significant} in announcing the subject at hand. Let me mention some of them, where the reference to Aristotle or Aristotelian ideas is obvious:

 \begin{itemize}

\item   \emph{Mathématique et réalité: faut-il croire Aristote ?} (Mathematics and reality: should we believe Aristotle?) \cite{realite}  (1986).

\item  \emph{Considérations sur la finalité: d'Aristote à Geoffroy Saint-Hilaire} (Considerations on finality: From Aristotle to Geoffroy Saint-Hilaire) \cite{Considerations} (1986).

 \item  \emph{Les intuitions topologiques primordiales de l'aristotélisme} (The primary topological intuitions of Aristotelism) \cite{Thom-Intuitions}  (1988).

\item Structure et fonction en biologie aristotélicienne (Structure and function of Aristotelian biology) \cite{Structure} (1988).

 \item \emph{Esquisse d'une sémiophysique: Physique aristotélicienne et théorie des catastrophes} (Semiophysics, a sketch: Aristotelian physics and catastrophe theory)
 \cite{Esquisse} (1988).

\item   \emph{Causality and finality in theoretical biology: a possible picture} \cite{Causality} (1989).

\item  \emph{Homéomères et anhoméomères en théorie biologique d'Aristote à aujourd'hui} (Homeomerous and anhomeomerous  since Aristotle until today)   \cite{Homeo} (1990).

\item  \emph{Matière, forme et catastrophes} (Matter, form and catastrophes)  \cite{A-Matiere} (1991).

\item \emph{Actualité de la physique aristotélicienne} (Actuality of Aristotelian physics) \cite{Actualite} (1992).

\item  \emph{Biologie aristotélicienne et topologie} (Aristotelian biology and topology) \cite{Biologie} (1995).

\item  \emph{Qualitative and quantitative in evolutionary theory with some thoughts on Aristotelian Biology}  \cite{Qualitative}  (1996).

\item \emph{Histoire de la transversalité d'Aristote à la théorie des
catastrophes} (History of transversality, from Aristotle to catastrophe theory)  \cite{Histoire} (1996).  

\item \emph{The hylemorphic schema in mathematics} \cite{hylemorphic}  (1997).

\item \emph{Comment la biologie moderne redécouvre la kinèsis d'Aristote} (How modern biology rediscovers Aristotle's \emph{kinesis})  \cite{Comment} (1998).

 \item \emph{Aristote topologue} (Aristotle as a topologist) \cite{Aristote-Topologue} (1999).

 \end{itemize}

   I will highlight several ideas that lie behind these works in the next section.\index{Aristotle}\index{Thom, René} Here I will only make a few comments on some papers in this list.
   
   In the article \emph{Actualité de la physique aristotelicienne} \cite{Actualite}, Thom starts by pointing out a few principles from Aristotelian physics\index{Aristotle!physics} that can be found in an equivalent form in catastrophe theory.\index{Thom, René}\index{catastrophe theory} He dwells on the distinction that Aristotle makes between the phenomena that take place in the celestial, or ``supra-lunar" world, and those of the sub-lunar one. The former is governed by God-given uniform circular motions, whereas the latter is subject to one great principle, namely, the effort of the four elements (fire, air, water and earth) which moves them toward what Aristotle calls their ``natural place".  It is in this distinction that Thom finds the origin of Aristotelian physics.\index{Aristotle!physics}    He makes a comparison between this division and the modern classification of dynamical systems into reversible Hamiltonian ones and irreversible gradient-like dissipative ones. This leads him to present Aristotelian physics in an axiomatic way. He starts with an axiom of existence of entities (\emph{ousiai}), continuing with their division into primary and secondary ones: The first ones include the living being, identified with a 3-dimensional topological ball separated from the exterior world by a boundary (the skin), and  the second ones are only ``qualities" that affect the primary ones. The other axioms concern the \emph{state} of an entity, the notion of \emph{power}, and the \emph{natural} transformations between entities. Entities interact with each other, and qualities are classified according to \emph{genus}.\index{Aristotle!genera} The analogy\index{analogy!René Thom} between the major concepts of  Aristotelian physics\index{Aristotle!physics} and those of catastrophe theory is outlined in terms of a few correspondences and principles:
 
 \begin{enumerate}
 
 \item Natural place $\longleftrightarrow$ genericity;
 
 \item genus decomposition  $\longleftrightarrow$ bifurcation;

 \item homeomeroous/anhomeomerous $\longleftrightarrow$
  regular and irregular point in a morphology;\index{morphology}
 
 \item biological organisation 
 $\longleftrightarrow$  stratification;

  \item The continuity assumption;
  \item Form is the boundary of matter.
 \end{enumerate}

  The fact that in Aristotle's \emph{Physics},\index{Aristotle!physics}  form is the boundary of matter, is an idea that Thom interprets using his background as a topologist. He writes \cite[p. 93]{Actualite} that the ambition of catastrophe theory is to keep this Aristotelian flexible scheme. We shall dwell on this in the next section.
  
   I will review now some ideas of Thom on form.

 Thom was, like  Aristotle, a philosopher of form.  Understanding form and its evolution, which Thom connected with his research in the topology of manifolds and function spaces, became an obsession for him in his philosophical works. He highlighted several passages in which his predecessor discusses form. 
In the \emph{Metaphysics}, Aristotle\index{Aristotle!metaphysics} writes that  what makes the essence of an object is its \emph{form}  \cite{A-Metaphysics} (1032b1-2, 1035b32). In the \emph{De Anima}, he  maintains that the soul is the form of a living being\index{Thom, René}  \cite{A-Soul} (412a11). 
In the \emph{Physics}, he writes that potentially, flesh or the bone do not have yet their own nature and do not exist by nature, as long as they have not received the form of flesh and bone, that is, a definable form,\index{form} that which can be stated to define the essence of the flesh or the bone \cite{A-Physics} (193b).

In the \emph{Esquisse}, Thom talks about Aristotle's theory of hylemorphism,\index{hylemorphism} to which he adheres completely. This is a theory which considers that any being (object or individual) is composed in an indissociable way of  matter (\emph{hylé}, \tg{>'ulh}) and of form\index{form} (\emph{morphè}, \tg{morf'h}), matter being,\index{matter} for him, essentially something in power, a substrate which awaits to receive a form, in order to become a substance, the substance of being, or being itself. 
Thom notes that the word \tg{>'ulh} was originally used to designate shapeless wood, and that it is to Aristotle himself that we owe its introduction in philosophy, while the latter was putting forward the notion of continuity, both for matter and for form.  Thom writes on this topic  \cite[p. 12]{Esquisse}: 
\begin{quote}\small
Of course, I knew that the hylemorphic\index{hylemorphism}\index{Thom, René} scheme which I use in the formalism of catastrophes had its origin in the work of the Stagirite. But I was unaware of the essential point, namely, that Aristotle had tried, in the \emph{Physics}, to build a theory of the world based not on the number, but on the continuum. He had thus realized (at least partially) the dream I have always had of developing a ``Mathematics of the Continuum" which takes the continuum as its starting point, if possible,  without any recourse to the intrinsic generativity of number.
  \end{quote}
   
    Talking again about form in Aristotle,\index{Aristotle} Thom writes in \cite[p. 53]{Predire}: 
  \begin{quote}\small     I see matter in an Aristotelian perspective, a kind of continuum that can acquire form. Form can be external, visible, or internal. The internal form is what one would call, from the semantic point of view, a quality. Aristotle's \emph{materia signata} is a matter provided with quality. I think that any quality can be seen precisely, to a certain extent, as a spatial form, an extended form in an abstract space. The original matter,\index{Aristotle!matter} if I may say so, is a bit like Aristotle's raw material, it is the substrate which can receive any kind of quality, any predicate.  The raw material is a kind of idealization which, very quickly, acquires qualities, forms.
 
         \end{quote}

Thom declared several times in his writings that he shared with Aristotle the conception of the soul as the form of the body. 
He comments on this in \cite[p. 100]{Predire}: 
\begin{quote}\small
The capacity of a body to be the support of a soul is something that presents itself as a structure or as a law of a formal nature, associated precisely with form in the morphological sense, in the spatio-temporal sense of organization. The whole is obviously linked to the form of the flows running through the organism: blood, neuronal, or more generally metabolic: for me all this is form, and something comes out of it as a residual form with an organizing character: the soul.  But the structure, the character (intrinsic in a certain sense) of this residual form, like the result of a formal structure, originates in this gigantic object that is the consideration of the form of all the molecular and physical movements of our organism.
\end{quote}

\section{Topology from an Aristotelian perspective}\label{s:topology}

 In his\index{topology} books and in several articles, Thom highlighted and  commented on several passages in the works of the Stagirite which support his conviction that the latter was aware of basic notions of topology. He declared in particular that a careful reading of Aristotle's \emph{Physics} shows\index{Thom, René} that the latter understood the topological distinction between a closed and an open set. The reader might try to imagine what kind of intuition can one have on the notion of open and closed set and the distinction between them if one does not have the language of set theory and topology.
 Thom writes in his \emph{Esquisse} \cite[p.\,167-168]{Esquisse}:         
    \begin{quote}\smaller 
          A careful reading of the \emph{Physica} leaves little doubt but that [Aristotle] had indeed perceived this difference. ``It is a whole and limited; not, however, by itself, but by something other than itself" (207a24-35)\footnote{Here and in the following,  while I am giving the relevant passages in Aristotle, I am translating Thom's interpretation of the Greek.} could hardly be interpreted except in terms of a bounded open set. In the same vein, the affirmation: ``The extremities of a body and of its envelope\index{envelopes} are the same" (211b12) can be identified, if the envelope is of a negligible thickness, with the well-known axiom of general topology: ``Closure of closure is closure itself" expressed by Kuratowski at the beginning of this century. This allows the Stagirite\index{Aristotle} to distinguish two infinites: the great infinite that envelops everything and the small infinite that is bounded. The latter is the infinite of the continuum, able to take an infinity of divisions (into parts that are themselves continuous). Whence the definition he proposes: ``The infinite has an intrinsic substrate, the sensible continuum" (208a).
            \end{quote}

Chapter 7 of Thom's \emph{Esquisse}\index{Thom, René} is titled  \emph{Perspectives in Aristotelian biology}. The first part concerns topology\ and bears the title \emph{The primordial topological intuitions of Aristotelianism: Aristotle and the continuum}. It starts as follows \cite[p.\,165]{Esquisse}: 
        \begin{quote}\small
        We shall present here those intuitions which we believe sub-tend all Aristotelianism. They are ideas that are never explicitly developed by the author, but which---to my mind---are the framework of the whole architecture of his system. We come across these ideas formulated ``incidentally" as they were condensed into a few small sentences that light up the whole corpus with their bright concision.
        \end{quote}

 In his 1988 article \emph{Les intuitions topologiques primordiales de l'aristotélisme} \cite{Thom-Intuitions} and in his 1991 article \emph{Matière, forme et catastrophes} \cite{A-Matiere}, Thom returns to these matters, explaining that the modern topological distinction between an open and a closed set is also expressed in Aristotle's distinction between form and matter, and between actuality and potentiality. We mentioned this concerning his distinction between actual and potential infinity. He declares that the difference between the notion of bounded and unbounded open set in Aristotle's philosophical system lies in the fact that the former may exist as a substrate of being whereas the latter cannot  \cite[p.\,396]{Thom-Intuitions}.
He notes the presence of the notion of boundary in formulae such as: ``Form\index{form} is the boundary\index{boundary} of matter,"\index{matter} \cite[p.\,398]{Thom-Intuitions} and ``Actuality is the boundary of potentiality" (\cite[p. 399]{Thom-Intuitions} and \cite[p.\,380]{A-Matiere}). He points out that the paradigmatic substance for Aristotle, which is the living being, is nothing else than a ball in Euclidean space whose boundary is a sphere, that is, a closed surface without boundary. Shapeless matter is enveloped by form---\emph{eidos}---in the same way as the boundary of a bronze statue defines its form.\index{form} The boundary of a living organism is its skin, and its ``interior" exists only as a potentiality. 

 The distinction between openness\index{open set!Aristotle} and closeness\index{closed set!Aristotle}\index{Aristotle!closed set}\index{Aristotle!open set}  in the topological sense in Aristotle's works is also pointed out by Thom in his discussion of homeomerous\index{homeomer} and anhomeomerous\index{anhomeomer} parts of animals, which we review in \S \ref{s:biology} below, in the context of Aristotle's biology.\index{Aristotle!biology}

 Another guiding idea in Thom's philosophical work which has an Aristotelian origin is the opposition between the discrete and the continuous.
He declares, about this opposition:  ``For me, the fundamental aporia of mathematics is indeed in the opposition discrete-continuous. And this aporia dominates at the same time all thought" \cite[p. 81]{Predire}.   In this interview, he recalls Zeno's paradox\footnote{I take this opportunity to point out an extremely interesting interpretation of Zeno's paradoxes by S. Negrepontis in \cite{Negrepontis}, based on a thorough interpretation of Plato's dialogue \emph{Parmenides}.\index{Parmenides} In short, Negrepontis' explanation is that Zeno's\index{Zeno} aim in his argumentation is to show that the sensible Beings are separate from the true/intelligible Beings.} supported by the story of Achilles and the tortoise and he declares:  
\begin{quote}\small
The continuous\index{continuity!philosophy of} is in a way the universal substratum of thought, and of mathematical thought in particular. But one cannot think anything effectively without having something like the discrete in this continuous unfolding of the mental processes: there are words, sentences, etc. The \emph{logos},\index{logos} the discourse, is always discrete; there are words entering in a certain succession, but discrete words. And the discrete immediately calls for the quantitative. There are points: we count them; there are words in a sentence: we can classify them quantitatively by the grammatical function they occupy in the sentence, but there is still an undeniable multiplicity.
\end{quote}

 Thom basically believed in the continuous\index{continuity!philosophy of} character of the universe,\index{Thom, René} of phenomena and of their substrate. Discontinuities only operate on the continuous. He was opposed to what he used to call a ``modern vulgate", which comes essentially from information science and which says that everything is expressed in bits (\cite [p. 62]{Predire}). In the \emph{Esquisse d'une sémiophysique} \cite[p. 172]{Esquisse}, there is a long passage in which he discusses Aristotle's decision to leave Plato's Academy, a decision which he explains by a disagreement with the master, concerning the notion of continuous. In particular, Thom writes that Aristotle's ``revolt against his master is that of the topologist against the imperialism of the arithmetician."
He writes:   ``Aristotle does not go into much detail about the way he sees the raw material. In my opinion, it must always be reduced to a continuum, to an extent. I am a universal topologist! I have a real metaphysics of the continuum"  \cite[p. 54]{Predire}.
           
           The notion of boundary\index{boundary!Aristotle}\index{Aristotle!boundary} is another key idea which Thom found in Aristotelian physics. He writes in \cite[p. 110]{Predire}:
  
\begin{quote}\small
 For Aristotle, the form\index{form!Aristotle}\index{Aristotle!form} of a physical object is something like its boundary; in the abstract sense of the definition, the \emph{eidos} is also something like a form in an abstract space, with its boundary. There is  intelligible matter, which in some way is compressed by its definition. \emph{Orismos} means definition; it is almost the same word as \emph{oros} which means boundary.\index{Thom, René} This is quite remarkable. \ldots To define is to draw the borders. I think that there's something quite profound about that intuition.
\end{quote}

 In the two articles 
 \emph{The primary topological intuitions of Aristotelism} \cite{Thom-Intuitions}  and 
\emph{Aristotle as a topologist} \cite{Aristote-Topologue}, Thom pushed\index{Thom, René} his ideas to the point of finding Stokes' famous\index{Stokes formula!Aristotle}\index{Aristotle!Stokes formula} formula in the work of the Stagirite.\index{Aristotle} He explained in particular how he found this formula in a notion of minimal limit of an ``enveloping body''  contained in Aristotle's Physics (211b 11-12), a formula that he translates as: ``The edge of the enveloped body and the edge of the enveloping body are together''. This identification of Stokes' formula with a formula of Aristotle is echoed again by Thom in his 1991 article \emph{Matière, forme et catastrophes} (Matter, Form and Catastrophe) \cite[p. 381]{A-Matiere}, an article published in the proceedings of a conference held at the Unesco headquarters in Paris at the occasion of the 23$^{\mathrm{rd}}$ centenary of the Philosopher, 
     If we add to this that Thom foresaw in the latter the mathematical ideas of genericity, stratification and bifurcation, we will not be surprised to see the reference to Aristotle in the subtitle of the \emph{Esquisse} \cite[p. 13]{Esquisse}.

The notion of cobordism,\index{cobordism} which was at the center of Thom's work in the 1950s, is also incorporated into this philosophical framework \emph{via} the notion of boundary, which is its primary element. He states in his 1992 interview \emph{La théorie des catastrophes} \cite{INA}: 
\begin{quote}\small
The whole unity of my work revolves around the notion of boundary, because the notion of cobordism is only a generalization of it. The notion of boundary seems all the more important to me today as I am interested in Aristotelian metaphysics. For Aristotle, the boundary is a principle of individuation. The marble statue is matter in the block from which the sculptor drew it, but it is its boundary which defines its form.
\end{quote}

        \section{Aristotelian biology}\label{s:biology}
 Like Aristotle, Thom was the philosopher of biology.\index{biology!Aristotle} And basically, biology, from Thom's point of view, was not far from topology, nor from philosophy: ``The great problem of biology," he says in a 1992 interview on catastrophe theory \cite{INA}, ``is the relationship between the local and the global. Knowing the organism, why the organs have this or that property at this point and not another, etc., that is the big problem". He then\index{Thom, René} adds:  
\begin{quote}\small
The problem of the relation between the local and the global is a problem that is essentially philosophical, that has to do with extent; this problem is the object of topology. Topology is essentially the study of the means to make the junction between a known local situation and a global property to be found, or vice versa: knowing a global property of space, find the local properties, around each point. There is a kind of deep methodological unity between topology and biology.
\end{quote}

  In his article \emph{Les mathématiques et l'intelligible} (The mathematics of the intelligible) (1975) \cite{T1975}, Thom introduces what he calls mathematical thought in the context of biology and anthropology,\index{Thom, René} a thought that, according to him, goes beyond the usual language in describing external phenomena with their mixture of determinism and indeterminism.  

  Chapter 5 of Thom's \emph{Esquisse} is titled \emph{The general plan of animal organization}. It is in the lineage of the zoological treatises of Aristotle, expressed in the  language of modern topology. Thom writes in the introduction:  
``This presentation might be called an essay in transcendental anatomy; by this I mean that animal organization will be considered here only from the topologists' abstract point of view." He adds: 
``We shall be concerned with ideal animals, stylized images of existing animals, leaving aside all considerations of quantitative size and biochemical composition, to retain only those inter-organic relations that have a topological and functional character."  In  \S B of the same chapter, Thom returns to Aristotle's notion of homeomerous and anhomeomerous, in relation with the stratification of an animal's organism, formulating in a modern topological language the condition for two organisms to have the same organisation. 
 Aristotle considers this notion at the beginning of the \emph{History of Animals} \cite{A-HA}: ``Regarding the parts of the animals, some are non-composite: these are those which can be divided into \emph{homeomers}, like flesh is divided into flesh; others are composed: these are those which can be divided into \emph{anhomeomers}; the hand, for example, cannot be divided into hands, nor the face into faces."
 Thom explains that a homeomerous part of an animal (\tg{'omoiomer'hs}) has generally a boundary structure constituted by anhomeomerous parts, which implies that the substrate of a homeomerous part is, topologically speaking, not closed. Mathematics, philosophy and biology are intermingled in this interpretation of a living organism, following Aristotle, with formulae like: 
``The opposition homeomerous-anhomeomerous is a `representation' (a homomorphic image) of the metaphysical opposition: potentiality-act. As the anhomeomerous  is part of the boundary of a homeomerous of one dimension higher, we recover a case of the application of \emph{act as boundary of the potentiality}."  \cite[p.\,400]{Thom-Intuitions}.

In Thom's words, an animal organism is a three-dimensional ball $O$ equipped with a stratification which is finite provided we neglect the details that are too fine. For example, when we consider the vascular system, we may take into account only the elements that can be seen with the naked eye: arterioles and veinlets. Seen from this point of view, the homeomerous\index{homeomer} parts are the three-dimensional strata: blood, flesh, the inside of bones, etc., the two-dimensional strata are the membranes: skin, mucous membrane, periosteum, intestinal wall, walls of the blood vessels, articulation surfaces, etc., the one-dimensional strata are the nerves:  vessel axes, hairs, etc., and the zero-dimensional strata are the points of junction between the one-dimensional strata or the punctual singularities: corners of the lips, ends of hairs, etc.

    This leads Thom\index{Thom, René} to the definition of two organisms $O$ and $O'$ ``to have the same organisation"  by the existence of a homeomorphism $h:O\to O'$  that preserves this stratification.

Thom also developed his ideas on the stratification of an animal body  in the series of lectures given in 1988 at the Solignac Abbey \cite{Essentielles}.
On p. 7 of these notes, he addresses the question of when two animals have isomorphic stratifications, and he uses for that the notion of isotopy between stratified spaces: Two sets $E_1, E_2$ have isotopic stratifications\index{stratification} if there exists a stratification of the product $E\times [0,1]$ such that the canonical projection $p:E\times [0,1]\to[0,1]$ is of rank one on every stratum of $E$, with $E_1=p^{-1}(1)$ and $E_2=p^{-1}(2)$. He considers that this notion is implicitly used by Aristotle in his classification of the animals, insisting on the fact that the latter neglected all the quantitative differences and was only interested in the qualitative ones.

Thom also talks about stratification in Aristotle in chapter 5 of his \emph{Esquisse}, a chapter titled \emph{The general plan of animal organization}, in relation with the notion of homeomerous\index{homeomer} introduced by the Stagirite.\index{Thom, René} 
He gives this idea of stratification in biology  an official status in his work,  by defining his \emph{phenomenological equivalence relation}, a relation which at the same time allows us to take a new look at Aristotle's classification of the elements of a living body. He writes in his article \emph{Structure et fonction en biologie aristotélicienne} \cite{Structure}, that this allows ``the identification of two anomalies by an isotopy between stratifications, either according to genus, or refined according to species."
  
In relation with Aristotle's classification of the parts of the animals, and more generally of a general plan of animal organisms, Thom revived a famous controversy that took place in 1830 between Geoffroy Saint-Hilaire and Cuvier, in which the former used a purely geometrical classification principle of animal organs whereas the latter defended criteria based only on function; see \cite{Thom-Archétype}, an article in which Thom proposes an approach to this dilemma using catastrophe theory.

  Genericity\index{genericity} is another idea that Thom finds in Aristotle's biology. He declares in his 1992 interview \emph{La théorie des catastrophes}  \cite{INA} that Aristotle considered that what matters in the study of nature is not to know everything that can happen, but everything that happens most often. In the Metaphysics, as in other texts, the Philosopher speaks at length about the difference between what happens most often and what is only an accident (\tg{sumbebhk'os}). He writes in particular that all  science is either about what is always or what is most often; otherwise, he says, how could one learn or teach others?\footnote{\emph{Metaphysics} \cite{A-Metaphysics} 1027a20.} Thom establishes a parallel between this and  the theory of singularities:\index{singularity} we study generic singularities, and those can be classified; their number is. There are important pathological singularities, but since they are not generic, the general theory does not deal with them.

 In the series of interviews \emph{To predict is not to explain},\index{predictive science} Thom returns to the homological Stokes formula $d^2=0$ and to its biological interpretation in terms of Aristotelian zoology: ``The boundary  of the boundary is empty; this is the great axiom of topology, of differential geometry in mathematics, but it expresses the spatial integrity of the edge of the organism." \cite[p. 111]{Predire}

Thom\index{Thom, René} spent some time trying to teach biologists the fundamental notions of topology that he thought they may need. Several sessions of the courses he gave at Solignac were devoted to this subject \cite{Essentielles}. His 1982 lecture \emph{Les réels et le calcul différentiel, ou la mathématique essentielle} \cite{essentiel} at the \'Ecole Normale contains a summary of what is needed for the biologist as a minimal mathematical background.

  In the next section, we talk more about the status of singularities in Thom's philosophical system.

  \section{Philosophy of singularities}\label{s:singularities}
  
  The goal of Thom's article\index{singularity!philosophy of} \emph{Philosophie de la singularité} (1988)\cite{Philosophy} is to show how the mathematical notion of singularity is used in our perception and description of the world. Not only this, but from Thom's point of view, our perception is essentially a perception of form, a transfer of topological structure from matter to our mind.
Thom starts by making a connection between the mathematical notion of analytic continuation in the holomorphic setting and prediction in science. In complex analytic geometry, analytic continuation allows, from the knowledge of a function on some open set, to extend it to a larger set, the domain of holomorphy\index{holomorphy} of the function. Thom views this as a sort of prediction: we predict the value of the function on a larger subset. He notes that what prevents a function from being extended further, depends on the singularities, on the boundary of the domain of holomorphy.

 Another well-known notable fact is that a  holomorphic function is determined by its singularities.\index{singularity!holomorphic function} This was already known to Bernhard Riemann. One form of the latter's famous \emph{Existence theorem} says that one can reconstruct a meromorphic function  from the knowledge of its singularities and its topological behavior at these singularities.  
   Paul Montel,\index{Montel, Paul} who has been Henri Cartan's doctoral advisor, emphasized the importance of characterizing a function of a complex variable in terms of its singularities. 
  In the introduction to his \emph{Leçons sur les fonctions enti\`eres ou méromorphes} \cite{Montel}, he writes:  ``What characterizes a function, what distinguishes it from others, is the set of its singularities.\index{Thom, René} What holds for functions also holds for  concrete objects and living beings" \cite[p. vi-vii]{Montel}.\footnote{Ce qui caractérise une fonction, ce qui la distingue des autres, c'est l'ensemble de ses singularités. Il en est des fonctions comme des objets concrets et des êtres vivants.}
  
        Another idea related to the fact that the singularities of a function give global information on it comes from Morse theory,\index{Morse theory} one of whose fundamental principles is that a Morsefunction allows the reconstruction of a space out of some finite number of cells of a fixed type and the singularities of the function. Such ideas were part of the background that  accompanied Thom throughout the rest of his life.

          It is during his time in Strasbourg that Thom began to reflect on Morse theory, which, in fact, is the subject of his first article, a note in the Comptes Rendus \cite{T-CRAS} published in  1949.\footnote{Haefliger, in his article \cite{H-IHES}, writes, before presenting the ideas contained in this note and its later developments: ``This was not too early, because the administrators at CNRS wondered if it was necessary to continue to support this young mathematician who was so little productive. [Il était temps, car les responsables du CNRS se demandaient s'il fallait continuer à soutenir ce jeune mathématicien si peu productif.]} What fascinated Thom in this theory was, once again, the fact that one can use the singular points of a function to decompose a manifold into cells that have a standard type, and conversely, one can reconstruct the space from the information on the cells and their order in the values of the function. Talking about this kind of situation, Thom recalls an old formula of Descartes:\index{Descartes, René} ``To reduce everything to situations that are simple enough to be described" \cite{INA}. This is another idea that never ceased to fascinate him: how to make constructions out of simple things. Regarding simplicity, he writes again in \emph{Paraboles et catastrophes} \cite[p. 66]{Paraboles}: ``The whole `philosophy' of catastrophe theory, its general scheme, is precisely this: \emph{ It is a hermeneutic theory which tries, in front of any experimental data, to build the simplest mathematical object that can generate it}."

        In the article \emph{Philosophie de la singularité}\index{singularity!philosophy} that we already mentioned, Thom motivates the use of analytic functions in physics, and the fact that they constitute the basis for prediction (via analytic continuation) by the fact that the symmetry groups of the laws of physics, such as those that govern particles and their interactions, are analytic. He also recalls that the notion of stratification, in mathematics, arises precisely in an analytic setting. Singularities govern the propagation of analyticity and at the same time they play the role of obstacle for this propagation.
        
      Regarding singularities, we also mention that   Thom also studied the notion of apparent contour of objects, contours which present singularities that he can classify. Observing apparent contours reminds him of the inhabitants of Plato's cave, chained and unable to move, 
  only capable of watching shadows of world events that take place outside the cave, projected on the wall in front of them  by the light of a fire burning at the entrance of the cave. He says that this is how we should interpret Hironaka's desingularization\index{Hironaka desingularization theorem} theorems.\footnote{In the 1970s, Heisuke Hironaka built a theory  on the desingularization of  of sub-analytic sets, showing that locally, any such set is a finite union of images of spheres  glued together by real-analytic maps. The theory had been foreseen by Thom.\index{Thom, René} See Hironaka's exposition in \cite{Hironaka}.} He notes incidentally that the singular sets of apparent contours are naturally stratified spaces in which a hierarchy of singularities is manifested, and he mentions a relation between such a situation and the notion of unfolding, a central notion in his theory of elementary catastrophes. Talking about the biological and human world, he mentions the ``archetypal singularities" that appear in embryology. He concludes the article by pointing out an opposition between singularities created as defects of an ambient propagative structure and singularities as a source of the propagative effect, a problem that, he says, is at the source of all scientific disciplines.

\section{Topology and ethology}\label{s:ethology}
In this short section, I briefly review the point of view of a researcher in bioethics, Yanick Farmer on the importance of Thom's morphological theories in epistemology and in the theory of ethics. This view is expressed in the article \emph{Topologie et modélisation chez René Thom : l'exemple d'un conflit de valeurs en éthique} (Topology and modelization after René Thom: The example of a conflict of values in ethics)\index{ethics} (1910)\cite{Farmer}.

Farmer considers that the abstract (topological) setting of Thom's approach to science, based uniquely on morphology,\index{morphology} answers the problem of the need for universality  in the particular setting of ethology,\index{ethology} by avoiding the problem of relativisation that is inherent to the use of specific languages that are proper to particular cultures. The author mentions as an example the notion of \emph{autonomy}, which is central in bioethics but which has different meanings according to the language and the culture where it is used.
Thus, from the ethologist point of view, to study a particular phenomenon, one must first differentiate it from its environment, and here, the author, like Thom, speaks of the need for this phenomenon or event to have a \emph{boundary},\index{boundary} and the \emph{substrate}\index{substratum} of the discussion becomes a mathematical space which may be assimilated to a 3-dimensional Euclidean space.

 I refer the interested reader to Farmer's article, in which he develops his ideas on the subject, which constitute a step towards  the  use of Thom's theory in ethology, but I would like to quote the conclusion of this article, where the author refers again to Plato's\index{Plato} cave passage (I am translating again from the French) \cite[p. 386]{Farmer}:
\begin{quote}\small
The philosophical positions that stand in the background of  Thom's theory of models seem to me to be faithful to the essence of the spirit of philosophy conceived by Plato\index{Plato} in his allegory of the cave. They act like a liberation for those among us who,  by the chains of knowledge and perceptual limitation, are bind to the surface of familiar objects described by ordinary language. For these prisoners who do not have at their disposal any other impressive technical mean than the strength of their gaze, Thom's philosophy offers an epistemological way of salvation. By challenging the methodological omnipotence of Reductionism, which forces knowledge to reveal the dephtless secret of elementary entities, René Thom reminds us, like Plato before him, that  the thinker of  everyday life and of the familiar things, as long as he is willing to become a geometer, can also try to understand better, in order to act better.
 
\end{quote}

    \section{Linguistics and morphology} \label{s:linguistics}
          For René Thom,\index{linguistics}\index{morphology} every science is  determined by a morphological scheme, and linguistics is no exception. In his hands, linguistics becomes a topological field. This is the subject of the present section.

      Thom was particularly interested in linguistics.  He saw in this field, as in biology, a privileged ground for the application of his ideas on morphology, with simple archetypes, where words are considered as living organisms, starting with the initial cry of the newborn, passing through the first stammering and stuttering, and culminating in the articulated and reasoned language of the adult. Aristotle, who studied language in his treatises on animals, considered that it is this articulated language, and the fact that it makes it possible for man to hold a reasoned discourse, that separates man from the other animals.\footnote{In the \emph{Politics} \cite{Politics} 1253a9-10, Aristotle considers that what separates man from other animals is the fact that man is a reasonable animal,  \tg{l'ogon >'eqon}, ``endowed with logos", the word \emph{logos} meaning here speech, or reasoned speech.} In this respect, Thom also spoke of language as a ``boundary", which, as we already pointed out, is a notion that accompanied him during all his life:  Following Aristotle,\index{Thom, René}\index{Aristotle} he considered language as the boundary that  separates man from the rest of the animals.

      In 1971, Thom published his article \emph{Topologie et linguistique} \cite{TL},  which also appeared five years later in Russia, translated by Yuri Manin. The article concerns the application of catastrophe theory to linguistics.
       Thom starts by exposing the idea developed by Ferdinand de Saussure, that linguistics is a system of signs that associates two elements: the ``signifier", which, here, is the sound, or the written word, and the ``signified", which is of psychological nature, and which is the mental representation of the concept associated with the sign. 
      The signifier\index{signifier} and the signified\index{signified} are considered by Thom as two morphological processes, and the whole approach of Thom is based on a transformation between morphologies, which allows us to give an explanation of the syntactic structure of a sentence thanks to its signified. He writes: ``For lack of sufficient mathematical knowledge, semanticians did not get the idea that the morphology\index{morphology}\index{Thom, René} of the signified could be supported by a multi-dimensional space. Therefore  they were reduced to very poor descriptions, inspired by the logical analysis of language."
  His first objective was then to describe the morphological structure of the signified. For this purpose, he used a known model, the neuro-physiological model of Zeeman, which leads to a geometrical model of the succession of ideas in the ``stream of consciousness" of the philosophy of introspection, with the existence of a structurally stable attractor.
 To summarize, Thom's\index{Thom, René} major contribution in linguistics\index{linguistics} was certainly the fact of bringing into this field the morphological models that he used in a field like biology. We refer the reader to Thom's article \emph{La linguistique, discipline morphologique exemplaire} (Linguistics, an ideal morphological discipline)  \cite{Thom1974}.

 Regarding the significant\index{significant}\index{insignificant} and the insignificant, Thom liked to repeat this fragment of Heraclitus, the philosopher of absolute motion: ``The master whose oracle is at Delphi does not say, does not hide, but he signifies",\footnote{Fragment B93, cf.  \cite{Heraclite-W}.} which Thom paraphrases by: ``Nature sends us signs that we have to interpret".  
 In this regard, Thom often recalls Aristotle, for whom Art imitates Nature and not the other way around.

Let me note also, talking about language, form and Aristotle, that the question of form in language is also present in the \emph{Problems} \cite
 {A-Problemes-A}, a compendium attributed to Aristotle,\footnote{Although, as in well known, there is an uncertainly concerning the authorship of the whole treatise, there is not doubt on the fact that it belongs to the Aristotelian school and that it was inspired, if not completely written, by Aristotle himself.} made up of questions, some of them simple and others that are real open research problems, encompassing all topics of science (mathematics, biology, physics, optics, psychology, astronomy, etc.). In problem 23 of section XI of the collection, it is said that ``the voice is air receiving a certain form",---a form which often changes, and which then vanishes.  Similarly, problem 57 states that the voice of a human being---who, among the animate beings, is distinguished  by the fact that he possesses the gift of speech---can take a multitude of nuances and forms.

 Between Aristotle and Thom, I wish to mention another mathematician-philosopher, Marin Mersenne, whose \emph{Harmonie universelle} (1636) contains a book titled \emph{Treatise of voice and of singing} \cite [p.~1-88]{Mersenne}. This book consists of an exposition of the physiology of the voice, both in humans and in animals, in the pure Aristotelian tradition.
 The statement of Proposition IX\footnote{Mersenne's \emph{Harmonie universelle} is organized into propositions and proofs, like in mathematical treatises.} of this book is the following: ``Voice is the matter of speech, and only man can speak.'' In the demonstration of this proposition, Mersenne develops the Aristotelian idea of matter and form in speech, writing in particular: ``We make use of the voice to form speech, as sculptors make use of wood and stones to make images; for images or statues are made by the different figures which are given to the matter of which they are made: and speech is a harmonic perspective, to which the voice serves as a picture to receive all kinds of images, since words are the images of the notions of the mind." Speech, he adds, ``is the form, ornament, and perfection of the voice, which can only be formed and figured into speech by man, as speech can only be formed into discourse by the mind."  
   
All this is close to what Thom said four centuries later.

     \section{Psychoanalysis}\label{s:psycho}

   In his \emph{Sketch of a semiophysics},\index{psychoanalysis} Thom writes that one of his main objectives is ``to clarify the analysis of the original psychical mechanisms of our species" \cite[p. 16]{Esquisse}.
   
   In 1990, Michèle Porte defended a doctoral dissertation in Paris whose title is \emph{Psychanalyse et sémiophysique : Études épistémologiques en métapsychologie et en dynamique qualitative} (Psychoanalysis\index{psychoanalysis}\index{semiophysics} and semiophysics: Epistemological studies in metapshychology and qualitative dynamics)  \cite{Porte1}. In this dissertation, the author uses Thom's catastrophe theory and his ideas on qualitative dynamics, and in particular, structural stability, to try to reformulate metapsychology. Her main thesis is that Sigmund Freud's\index{Freud, Sigmund} theories cannot be understood in a formalist and structural setting, but, on the contrary, in a dynamical setting of birth, persistance and vanishing of forms and the dynamics behind them. Thus, she proposes an interpretation of Freud's ideas using a morphodynamical\index{morphodynamics} model of psychical events, in the tradition of Thom.\index{Thom, René} At the same time, she makes a confrontation between Freud's psychoanalysis and  Thom's semiophysics.\index{semiophysics}
   
 Later, the same author published a systematic essay  titled \emph{La dynamique qualitative en psychanalyse} (Qualitative dynamics in  psychoanalysis) (1994), of which  Thom wrote the preface. She writes \cite[p. 20]{Porte2}:
 \begin{quote}\small
 It is not forbidden to think that in the no man's land that our culture has constituted between the two slopes of the logos, that of number and that of word in general, besides catastrophe theory\index{catastrophe theory} and semiophysics,\index{semioophysics} a transposed psychoanalysis has its share and its role to play \ldots Thus, we give a reason for the ``exterritoriality" of psychoanalysis, a theme as hackneyed as it is obscure. ``When we have a thought, the meaning of that thought is the form of the underlying neurophysiological process", according to one of the formulae by which Thom has often paraphrased Bernhard Riemann.\index{Riemann, Bernhard} Thom has taken this statement seriously; Freud's\index{Freud, Sigmund} work does not object to it; we propose that both analysts and other researchers take up the program of work that this statement imposes.
 \end{quote}

 Talking about Riemann, let me mention that Thom was also a reader of the latter's philosophical fragments. In his article \cite{Tenerife}, he  writes: 
  \begin{quote}\small
  As for myself, it is in reading in Riemann's\index{Riemann, Bernhard} \emph{Philosophisches Nachlass} the sentence ``When we think about thought, the meaning of this thought is nothing but the `Form' of the underlying physiological process" that I see how the concept of Form can help bridge the gap between the two types of science [Geisteswissenschaften and Naturwissenschaften]. In this light we may reasonably evoke a revival of \emph{Naturphilosophie}.
\end{quote}

 Porte's reading of Freud using Thom's background brings a fresh point of view on the work of the Viennese Master.
 
Let me also mention that Thom was aware of the applications of catastrophe theory to psychoanalysis since the beginning of catastrophe theory. We already quoted his comments on Zeeman's\index{Zeeman, Christopher} article  \emph{Topology of the brain} \cite{Zeeman} which gave him an impetus for his ideas on biological modeling; see also Thom's comments in \cite{Thom-Logos} (1989). Zeeman's work included the application of catastrophe theory to psychology, ethology and sociology. In the same article, Thom writes that he noticed that in Riemann's  \emph{Complete works} edition, there are a few pages where the latter expresses his view on the relation between spirit and matter, and where he claims that one can conceive that if an activity of the brain  corresponds to an idea that we think, then the meaning is necessarily defined by the \emph{form} of the corresponding physico-chemical process. 
 
 Let me mention, to close this section, that Thom gave a talk at a conference on psychoanalysis titled ``L'inconscient et la Science", in which he outlined a possible use in psychoanalysis\index{Thom, René}\index{psychoanalysis} of his ideas on morphology in the way he applied it in biology and linguistics,\index{morphology}\index{biology}\index{linguistics} see \cite{Thom1991}.
 
 \section{A metaphysics of analogy}\label{s:analogy}
 
In the 1970s,\index{analogy}\index{analogy!René Thom} Thom started a systematic use of analogies and metaphors in his writings. In fact, his approach to phenomena of science and morphology as taking place in an abstract substrate favored an extensive use of analogy.
He declares in the interview \emph{To predict is not to explain}  \cite[p. 122]{Predire}: ``There could be science only if we immerse the real in a controlled virtual. And it is through the extension of the real into a larger virtual that one can study the constraints that define the propagation of the real within that virtual."\footnote{Il n'y a de science que dans la mesure où l'on plonge le réel dans un virtuel contrôlé. Et c'est par l'extension du réel dans un virtuel plus grand que l'on étudie ensuite les contraintes qui définissent la propagation du réel au sein de ce virtuel.} On a related subject,\index{analogy!René Thom} 
 talking about the notion of ``boundary'', he declares in another interview: ``To reach the limits of what is possible, you have to dream the impossible, and it is really the interface between the possible and the impossible that is important because if we know it, we know exactly the limits of our power"  \cite{Nimier}.\footnote{Pour atteindre les limites du possible, il faut rêver l'impossible, et c'est réellement l'interface entre le possible et l'impossible qui est important parce que si nous le connaissons, nous connaissons exactement les limites de notre pouvoir.}
 
 In the field of analogy,\index{analogy} Thom considered himself once again the heir of Aristotle. In a 1976 conference, titled \emph{Le statut épistémologique de  la  théorie  des catastrophes} \cite{A-statut}, he writes that  ``the simple fact of being able to classify analogous situations is a considerable  
  achievement: there has been no theory of analogy since Aristotle."  
 Later, in the interview titled \emph{La théorie des catastrophes} (1992) \cite{INA}, he declares that catastrophe theory\index{analogy!René Thom} is  ``a theory of analogy",\index{catastrophe theory}\index{analogy} the first, of such breadth, since Aristotelian logic.

From a strictly mathematical point of view, the word\index{analogy} ``analogia"\index{analogy}\index{Thom, René} is used by Euclid\index{Euclid} in the \emph{Elements} to denote an equality between two ratios. Thom quotes at several occasions  Aristotle's example\index{Aristotle!analogy}\index{analogy}\index{analogy!Aristotle} of analogy saying that old age is to life the same thing as evening is to day,  in the form of an equality between two ratios:
  \[\frac{\mathrm{old \ age}}{\mathrm{life}}=\frac{\mathrm{evening}}{\mathrm{day}}
 .\]
From this follows the metaphor:\index{metaphor} ``old age is the evening of life".\footnote{The example is often attributed to Aristotle,\index{Aristotle} but the latter, in a passage of the \emph{Poetics}, quoting this sentence, refers to Empedocles,\index{Empedocles} who said that old age is the sunset of life.} In the same passage, Aristotle writes: ``The cup is to Bacchus what the shield is to Mars. We shall therefore say: ``the shield, cup of Mars", and ``the cup, shield of Bacchus".\footnote{Cf. \emph{Poetics} \cite{Poetics} 1457b20-21, and \emph{Rhetorics} \cite{Rhetorics} 1412b35 and 1413a1-2.}

  In his \emph{Esquisse d'une sémiophysique} and in other works,\index{analogy} Thom  highlighted Aristotle's analogies concerning animals.\index{analogy!René Thom}\index{analogy!Aristotle}\index{Aristotle!analogy} Let us especially mention an analogy that Aristotle makes in his \emph{Parts of the Animals} (668a10-13), comparing the vascular system of an animal to an irrigation channel in a garden.  Thom  emphasized this analogy  in his 1988 lectures titled \emph{Structure and Function in Aristotelian Biology} \cite{Structure}.\footnote{A similar\index{analogy}  agricultural metaphor\index{Aristotle!analogy} was made by Plato\index{Plato} in the \emph{Timaeus} \cite{Timaeus} (77c7-8), about blood irrigation in humans.}  Among the other examples of analogies made by Aristotle are the one between the organization of the parts of an animal and the political organization of the city (\emph{Politics}, 1290b25-35), the one saying that sight is to the body what intelligence is to the soul (\emph{Nicomachean Ethics},  1096b28), the one saying that the soul is to the body what the pilot is to the ship (\emph{De Anima}, 413 a 8-9), the one saying that the scale is to the fish what the feather is to the bird (\emph{History of Animals}, 486b15), and there are many others.   
 Thom,\index{Thom, René} on several occasions,\index{analogy} spoke of an analogy\index{analogy!René Thom} between the development of the embryo and that of a Taylor series function.  I have commented on other instances of analogy as a method of thinking in the works of Aristotle and Thom in my article  \cite{A-Logos-et-analogie}.

Thom also talked about Aristotle's theory of analogy\index{Aristotle!analogy}\index{analogy}\index{analogy!Aristotle} in terms of boundary. He writes in \cite[p. 75]{Predire}:
\begin{quote}\small
 Analogy\index{analogy} is not something arbitrary. Analogy, or metaphor,\index{metaphor} contrary to the common vision that makes it something approximate, fuzzy, appears to me as a strict relation that,  in many cases, we can express mathematically, even if this mathematical expression in itself is not interesting in the mental process that makes you consider the analogy.
 \end{quote}
In the same passage, speaking again of Aristotle's\index{Aristotle!analogy}\index{analogy}\index{analogy!Aristotle} analogy between evening with respect to day, and old age with respect to life, he writes: ``The formal structure of this analogy is simply the notion of boundary. You have here a time interval; this interval has an end; one calls `evening' or `old age' the tubular neighborhood, if I dare say of the word `end', and the corresponding catastrophe is,  for me, a fold."

  Concerning the general role of analogy in science, Thom\index{Thom, René} writes, in his article \emph{The twofold way of catastrophe theory} \cite{L-Twofold}  that this role is very much contested and that it is the object of a vast misunderstanding among professional scientists, most of whom---imbued with a positivist spirit---see it as a source of illusion, claiming that science must be quantitative whereas analogy is a qualitative thought. He reminds the reader that Konrad Lorenz,\index{Lorenz, Konrad!analogy}\index{analogy!Konrad Lorenz} the famous biologist, in his reception speech of the Nobel Prize, stated that in a certain sense, all analogies are true, by definition. Thom, more precisely,\index{analogy!René Thom}  distinguishes two cases of analogies: a ``mathematically formalized analogy, associated with an organizing catastrophe", which is always ``true, but can be sterile (it can only generate more or less poetic metaphors)", and a ``partial" analogy, whose algebraization is incomplete, and he says about the latter: ``It is precisely by striving to specify an analogy that one can come across interesting data; it is the incompleteness of the analogy that offers the best possibilities of synthesis. It is only because we accept the risk of error that we can reap new discoveries."

                       Concerning analogy,\index{analogy}  I would like also to mention Kurt Gödel,\index{Gödel, Kurt} whom Thom often quoted for having shown that a formal system of arithmetic is either inconsistent or incomplete. Gödel made an analogy between this formal side of mathematics and a rigid human society. In a draft letter to David Scurlock, quoted by Pierre Cassou-Noguès in his book \emph{Les Démons de G\"odel : Logique et folie}  (Gödel's Demons: Logic and Madness) \cite[p. 164]{Cassou},\footnote{P. Cassou-Noguès quotes this letter as being part of the Gödel collection held in the library of the Institute of Advanced Study in Princeton.  The letter is also quoted by Hao Wang in \cite[p. 4]{Wang}.} Gödel writes that one may expect from a society without any freedom at all (i.e., proceeding in everything according to strict rules of  ``conformity") to be either inconsistent or incomplete in its behavior, i.e., unable to solve certain problems, perhaps of vital importance. He also says that a similar remark can be applied to human beings individually.  
 
   Talking about metaphor,\index{metaphor} I cannot but point out the collection of articles \emph{Mathematics as metaphor} by Manin \cite{Manin}.\index{metaphor}

  \section{Another return to the classics}\label{s:return}
      
Thom had the particular wish to make  topology and  the theory of singularities a universal framework which would provide models for all phenomena in the natural and social sciences. Was he convincing? That is another question, but it is a well-known fact that history of knowledge is full of examples of sketches of fundamental theories that took centuries to convince, because they were ahead of their time. 

  Among the mathematicians who preceded Thom and whose work had a strong philosophical dimension, one can think of Descartes, Leibniz, Euler, Lambert, Riemann, Poincaré, Grassmann,  Brouwer and Weyl, to mention only some of the most famous. All of them considered that the language of mathematics, and in particular that of geometry, goes beyond the mathematical sciences.  To mention, while remaining close to Thom's ideas, a famous case of a preeminent mathematician who proposed an idea that eventually was not understood by his contemporaries, one may think of Leibniz's\index{Leibniz, Gottfried Wilhelm}\index{universal characteristic} \emph{universal characteristic},\footnote{Leibniz used several names for the new field he had in mind:\index{universal characteristic}  \emph{analysis situs, geometria
situs, characteristica situs, characteristica geometrica, analysis geometrica, speciosa
situs}, etc.} a theory which he intended to develop in reaction to the methods of algebra, which did not satisfy him.  The universal characteristic is closely connected with topology in many ways. Leibniz described this field as being more general than that of rigid Euclidean geometry and of  Descartes's analytic geometry. He declared that this domain  would be purely qualitative and that it was concerned with the study of figures
independently of their metrical properties.

 Thom\index{Thom, René} refers to Leibniz's\index{Leibniz, Gottfried Wilhelm}\index{universal characteristic}  universal characteristic  in the book \emph{Morphogénèse et imaginaire} \cite [p. 23]{Morph}:
\begin{quote}\small
In fact, the ultimate ambition of catastrophe\index{catastrophe theory}  theory  is to abolish the distinction between mathematical language and natural language which has been rampant in science since the Galilean break. [\ldots] A geometrical modeling of ordinary verbal thought will only be of interest if we can, thanks to it, arrive at assertions that the usual logic of natural language cannot provide. This supposes that we are able to:
 
1) Modelize geometrically all the (rigorous) deductions of ordinary thought. In other words: to realize the Leibnizian dream of ``universal characteristic";

2) Go beyond this.
\end{quote}

Leibniz\index{Leibniz, Gottfried Wilhelm}  had brought together a number of ideas on this subject and submitted them to Christiaan Huygens, whom he considered as his mentor, but the latter did not find these ideas worth exploring.\footnote{The two men had met for the first time in 1672, in Paris, where Huygens was settled since 1666. Huygens was 17 years older than Leibniz, and for some time he was his private teacher in mathematics.} 
  Leibniz, who had an infinite admiration for his former master, abandoned the project.
   It is interesting to read here an excerpt from a letter written by Leibniz to Huygens, dated September 8, 1679: 

 \begin{quote}\small
 I am still not happy with Algebra, since
it provides neither the shortest ways nor the most beautiful constructions of Geometry.
This is why when it comes to that, I think that we need another analysis which is properly
geometric or linear, which expresses to us directly \emph{situm}, in the same way as algebra expresses \emph{magnitudinem}. And I think that I have the tools for that, and that we might represent figures and even engines and motion in words, in the same way as algebra represents numbers in magnitude. I am  sending you an essay, which seems to me worth considering; there is no other person than you, Sir,  who can better judge it, and your opinion will count  for me more than many others.
\end{quote}

Leibniz\index{Leibniz, Gottfried Wilhelm}  sent to Huygens\index{Huygens, Christiaan}  his manuscript titled \emph{Characteristica
Geometrica}, dated August 10, 1679.  The piece is published in Volume II of Leibniz' \emph{Mathematische Schriften}, edited by Gerhardt \cite[Vol. II, p. 141-168]{Leibniz-Gerhardt}.
and it is also reproduced in Vol. XIII of Huygens' \emph{Collected Works}  \cite[p. 219-224]{Huygens-8}. It starts with the words:

\begin{quote}\small

I found some elements of a new characteristic, completely different from Algebra and which
will have great advantages for the mental representation in an exact and natural manner, although without
figures, of everything that depends on the imagination. Algebra is nothing but the characteristic
of undetermined numbers or magnitudes. But it does not directly express place,
angles and motions, from which it follows that it is often difficult to reduce, in a computation,
what is in a figure, and that it is even more difficult to find geometrical proofs and
constructions which are enough practical even when the Algebraic calculations are done. But this new characteristic, based on figures we can see, cannot fail to give at the same time the solution and the construction of geometric proof, in a natural manner, and by analysis, that is, by determined ways.
Algebra is obliged to suppose the Elements of Geometry; instead, this characteristic pushes the analysis to the end: if it were completed in the way I conceive it, one could make in characters which will be only letters of the alphabet the description of a machine whatever compound it could be, which would give the means to the mind to know it distinctly and easily with all the parts and even with their use and movement, without using figures or models and without hampering the imagination.
 
\end{quote}

Leibniz\index{Leibniz, Gottfried Wilhelm}  then explains in more detail his vision of this new domain of mathematics,
and where it stands with respect to algebra and geometry, giving several examples of
a way to denote loci, showing how this allows the expression of statements such as the fact
that the intersection of two spherical surfaces is a circle, and the intersection of two
planes is a line.
The piece ends with the words \cite[p. 224]{Huygens-8}:

 \begin{quote}\small 
 I have only one remark to add, namely, that I see that it is possible to extend the characteristic
to things which are not subject to imagination. But this is too important and it would lead us
too far for me to be able to explain that in a few words. 

\end{quote}

  Huygens was not convinced by Leibniz's project. He  responded to him in a letter dated
November 22, 1679 ([176] p. 577):

 \begin{quote}\small 
 
 I have examined carefully what you are asking me regarding your new characteristic, and
to be frank with you, I cannot not conceive the fact that you have too much expectations
from what you spread on me. Your example of places concerns only realities that 
were already perfectly known, and the proposition saying that the intersection of a plane
and a spherical surface makes the circumference of a circle does not follow clearly. Finally,
I cannot see in what way you can apply your characteristic to which you seem you want
to reduce all these different matters, like the quadratures, the invention of curves by the
properties of tangents, the irrational roots of equations, Diophantus' problems, the shortest
and the most beautiful constructions of the geometric problems. And what still appears to
me stranger than anything else is the invention and the explanation of machines. I say it to
you unsuspiciously: in my opinion this is only wishful thinking, and I need other proofs in
order to believe that there could be some reality in what you present. I would nevertheless
restrain myself from saying that you are mistaken, knowing the subtlety and the deepness
of your mind. I only beg you that the magnificence of the things you are searching won't let
you postpone from giving us those which you already found, like this Arithmetic Quadrature
you discovered, concerning the roots of the equations beyond the cubical, if you are still
satisfied with it.\index{universal characteristic}\index{Thom, René} 
 \end{quote}

On continuity\index{continuity!philosophy of} in nature, we already mentioned that Thom constantly referred to Aristotle.\index{continuity!Aristotle}\index{Aristotle!continuity} Talking about Leibniz, we may recall that he declared \cite{Leibniz}: ``Nothing is done all at once, and it is one of my great and most verified maxims that nature never makes jumps. I called this the law of continuity, when I spoke of it in the news of the republic of letters; and the use of this law is very considerable in Physics,"\index{continuity!Leibniz}\index{Leibniz, Gottfried Wilhelm!continuity} a sentence that reminds us of what Thom declared three centuries later.

 The essay that Leibniz sent to Huygens remained practically unknown. But it drew the attention of a few 19$^{\mathrm{th}}$-century mathematicians,
including Hermann Grassmann\index{Grassmann, Hermann}  (1809-1877), who is considered as the founder of abstract multi-dimensional linear algebra, and, most of all, of a geometric interpretation of this algebra.
He was among the first to stress on the importance of Leibniz' \emph{Analysis situs}. We refer the interested reader to the articles \cite{Heath} by Heath and \cite{Eche} by Echeverr\`\i a.

   There exist two recent commented editions
of Leibniz's work on the geometric characteristic,\index{universal characteristic}  both included in doctoral dissertations, the one by J. Echeverr\`\i a \cite{Echeverria1995} (1995),
in France, and the other one by de Risi, \cite{Risi} (2007), in Germany.

 Like Thom did after him, Grassmann\index{Grassmann, Hermann} also stressed on the fact that mathematics is the adequate language for understanding the physical world. He formulated a ``universal geometric calculus"  which was successfully used in the foundations of mechanics, relativity, electrodynamics and quantum mechanics. 
  We refer the reader to the survey article 
 \cite{Hestenes} by D. Hestenes, in which the author gives a glimpse of Grassmann's vision and concludes by the words   ``An adequate history of geometric calculus remains to be written".
  
 We note in passing that Thom also talked about a universal grammar, see his paper \cite{Grammaire}.

\section{Conclusion}\label{s:conclusion}

Thom\index{Thom, René} was a singular figure in the history of mathematics, and one of his important contributions was to challenge all dogmas, especially those of science.
 The fact that he had shifted his focus to philosophy and the consequences of this change were not understood by most of his colleagues: criticisms were voiced here and there, by mathematicians and others, even though the Fields medal assured him an immunity from any dispute of his mathematical talents.  
         
It is not surprising that some colleagues of Thom considered his philosophical reflections as delirium. Some of the criticisms that came from the mathematical community were deeply unfair.
In the introduction to his book  on catastrophe theory \cite{Arnold}, Vladimir Arnold writes (p. {\sc viii}):  
\begin{quote}\small
Neither\index{Thom, René} in 1965 nor later was I ever able to understand a word of Thom's own talks on catastrophes. He once described them to me (in French?) as ``bla-bla-bla", when I asked him, in the early seventies, whether he had proved his announcements. Even today, I don't know whether Thom's statements on the topological classification of bifurcations in gradient dynamical systems depending on four parameters are true. [\ldots] Nor am I able to discuss other, more philosophical or poetical declarations by Thom, formulated so as to make it impossible to decide whether they are true or false (in a manner typical rather of medieval science before Descartes or (the) Bacon(s)).\footnote{It is fair to add that several years later, Arnold changed his mind concerning Thom. In an interview published in  avril 1997 \cite{Arnold-Notices}, he writes: 
``I am deeply indebted to Thom, whose singularity
seminar at the Institut des Hautes \'Etudes
Scientifiques, which I frequented throughout
the year 1965, profoundly changed my mathematical
universe. I was always delighted by the
way in which Thom discussed mathematics,
using sentences obviously having no strict logical
meaning at all. While I was never able to completely
free myself from the straitjacket of logic,
I was forever poisoned by the dream of the irresponsible
mathematical speculation with no
exact meaning."}
\end{quote}
       Thom did not waver, but deep down, the reactions to his ideas hurt him. He wrote in 
  \emph{La science malgré tout} (Science despite everything) \cite{EU} (1965): 
   \begin{quote}\small
   It used to be fashionable---and probably still is---in scientific circles, to rave about philosophy. And yet, who could deny that the only really important problems for man are philosophical problems? But philosophical problems, being the most important ones, are also the most difficult ones; in this field, to show originality is very difficult, a fortiori to discover a new truth. This is why the society, very wisely, has given up subsidizing research on philosophical subjects, where the profit is too uncertain, to devote its efforts to scientific research, where, thank God, it is not necessary to be a genius to do a ``useful work". 
      \end{quote}

Did the community of philosophers understand Thom?\index{Thom, René} In some sense, the answer is simple: unlike the situation in mathematics, there is no philosophical community.     In his article \emph{Leaving mathematics for philosophy}, \cite{Leaving}, Thom explains that in philosophy, things are not as simple as in mathematics: there are schools, cliques, etc.

  Regardless of his involvement in philosophy, Thom had a personal view on mathematics and on mathematicians and why we do mathematics, which he expressed sometimes in a very metaphorical\index{metaphor} mode. In his paper \emph{De l'icône au symbole. Esquisse d'une théorie du symbolisme} (From icon to symbol. A sketch of a theory of symbolism)\cite{Icone}, he writes:
   \begin{quote}\small

   I would like to see the mathematician as a perpetual newborn babbling in front of nature; only those who know how to listen to Mother Nature's answer will later manage to open a dialogue with her, and to master a new language. The others only babble\index{Thom, René}, buzzing in the void --- \emph{bombinans in vacuo}. And where, you may ask, could the mathematician hear the answer of nature? The voice of reality is in the sense of the symbol.
 
   \end{quote}

 Manin\index{Manin, Yuri} is among the mathematicians who understood Thom's approach since the beginning.  In his contribution to Kyoto's ICM  (1990), he writes \cite{Manin-ICM}:
 \begin{quote} \small [\ldots]   What is relevant is the imbalance between various basic values which is produced by the emphasis on proof. Proof itself is a derivative of the notion of ``truth." There are a lot of values besides truth, among them ``activities," ``beauty" and ``understanding," which are essential in the high school teaching and later. Neglecting precisely these values, a teacher (or a university professor) tragically fails. Unfortunately, this also is not universally recognized. A sociological analysis of the controversies around Catastrophe Theory of René Thom shows that exactly the shift of orientation from formal proof to understanding provoked such a sharp criticism. But of course, Catastrophe Theory\index{catastrophe theory}
    is one of the developed mathematical metaphors\index{metaphor} and should only be judged as such.
    \end{quote}

In any case, Thom\index{Thom, René} felt out of step with the mathematical milieu of his time, which was rather very monolithic. He writes \cite[p. 29]{Paraboles}:  
\begin{quote}\small
I never really thought of myself as a mathematician. In fact, a mathematician must have, in my opinion, a taste for difficulty, for beautiful, rich and deep structures. I do not have that taste at all. The ultra-refined structures that fascinate my colleagues---Lie groups, simple finite groups, etc., in short, all these kinds of mathematical mythologies--- have never really interested me. On the other hand, I like things that move, flexible things that I can transform to my liking. I prefer the field of mathematics where one does not really know what he is doing! This is the reason why I consider mathematics today with a certain detachment, and I cannot say whether there is currently a strictly mathematical problem for which I have a deep interest.
\end{quote}  

Thom\index{Thom, René} declared several times that  even though he became a mathematician  by accident, mathematics brought him great satisfaction. In the interview \cite{Predire}, he says:  ``Really, when one has found a theorem in his lifetime, he can tell that he participated in a certain form of immortality, whatever he did. An illusion perhaps, but among all the fictitious immortalities by which we are deluded, this is still one of the most solid" \cite[p. 73]{Predire}.   In the article \emph{La science malgré tout}\cite{EU}, he writes: ``I know only one  really difficult science: mathematics.''    

\vfill\eject

      \centerline{\sc Appendix: Jean Cavaillès and Albert Lautman}
      \bigskip

 Let me say a few words on the two philosophers of science, Jean Cavaillès\index{Cavaillès, Jean} (1903-1944) and Albert Lautman  (1908-1944) to whom Thom was attracted as a young student.

  Cavaillès is the author of an important philosophical corpus on the foundations of mathematics in which the stress is on a dynamical evolution of these foundations. Among the multitude of philosophical schools of the first half of the 20$^{\mathrm{th}}$ century, Cavaillès was situated between Hilbert's formalism and Brouwer's intuitivism, a trend that he used to call ``modified formalism". He was close to the mathematician \'Emile Borel, one of the French representatives of the ``semi-intuitionist" current. As a matter of fact, in 1940, Cavaillès wrote an authoritative synthesis of Borel's theory of the quantification of chance \cite{Cavailles-Borel}, which constitutes the latter's point of view on the philosophy of probability theory expressed in Vol. IV of his 
\emph{Traité du calcul des probabilités et de ses applications} \cite{Borel}. The volume is titled \emph{Valeur pratique et philosophie des probabilités} \cite{Borel}.  

Cavaillès taught at the \'Ecole Normale Supérieure before Thom entered there. His philosophical work, which was in part the result of a close collaboration with Emmy Neother, had a non-negligible impact on the mathematical research done at the \'Ecole.

Lautman had been Cavaillès's student at the \'Ecole Normale.\index{Lautman, Albert} 
In 1937, he defended a doctoral dissertation, in two parts,\footnote{The French doctoral system (doctoral d'état) required the defense of two dissertations. The first one was the main thesis, and the subject of the second one was proposed by the jury of the doctorate, a few months before the defense. Its preparation was supposed to take only a few months (usually 3 to 4).} the first one titled \emph{Essai sur les notions de structure et d'existence en mathématiques} (Essay on the notions of structure and existence in mathematics), and the second one,  \emph{Essai sur l'unité des sciences mathématiques dans leur développement actuel} (Essay on the unity of the mathematical sciences in their present development). Following the path of Cavaillès (and before him, that of Poincaré), Lautman considered that both the formalist and the intuitionist  movements were a failure. He was an advocate of structuralism in the tradition of Bourbaki.\index{Bourbaki, Nicolas} He was a promoter of the concept of unity of mathematics, see the collection of articles \cite{Lautman1}.

Cavailles and Lautman found the sources of their theories in the recent developments in mathematics and physics (notably quantum physics). They both wondered about the role and the consequences of the various movements of thought that had appeared at the end of the 19$^{\mathrm{th}}$ century on the philosophy of mathematics (conventionalism, logicism, constructivism, formalism, etc.). Even if they diverged on some points, and in particular on the organization of mathematics as a system of thought, they both considered that the philosophy of mathematics must necessarily be at the center of any metaphysical theory, as it was already  for Plato, Heidegger and other philosophers. Their work is an embodiment of this approach. Both Cavailles and Lautman had close relations with mathematicians of the Bourbaki group such as Cartan, Chevalley, Dieudonné, Ehresmann, and Weil and they were very much interested in Boubaki's\index{Bourbaki, Nicolas} project of writing complete treatises on the foundations of several topics in mathematics.

 During the Second World War,  Lautman et Cavaillès joined the resistance to the German occupation of France. They were both shot by the Nazis in 1944. The first was 36 years old and the second 40.

The French Society of Philosophy devoted its session of February 4, 1939 to the discussion of Cavaillès and Lautman's works. Several mathematicians, including Henri Cartan, Paul Dubreuil, Paul Lévy, Maurice Fréchet, Charles Ehresmann and Claude Chabauty, were present at that session. A report on this session was published after the liberation of France, in the Bulletin de la Société française de philosophie \cite{SFP} (1946). 
Thom may have been present at this debate, but being still too young, his name is not mentioned in the Annals.

  \printindex

     \end{document}